\documentclass[%
 reprint,
 amsmath,amssymb,
 aps,
nolongbibliography
]{revtex4}

\usepackage{graphicx}
\usepackage{dcolumn}
\usepackage{bm}

\usepackage{lipsum}
\usepackage{dsfont,eucal,mathrsfs}
\usepackage{verbatim}
\usepackage{color} 
\usepackage{float}
\usepackage{rotating}
\usepackage{url}
\usepackage{subcaption}
\usepackage{rotating}
\usepackage{adjustbox}
\usepackage{epsfig,amsmath,amsfonts,amsthm,graphicx}
\usepackage{bibunits}

\usepackage{amsmath}

\newcommand{\bo}{\boldsymbol}

\renewcommand{\mathcal}{\mathscr}

\newenvironment{bsmallmatrix}
  {\left[\begin{smallmatrix}}
  {\end{smallmatrix}\right]}
\usepackage{ragged2e}

\begin{document}
\begin{bibunit}
\title{Reconstructing Network Dynamics of Coupled Discrete Chaotic Units from Data}

\author{Irem Topal}
\email{irem.topal@khas.edu.tr}
\author{Deniz Eroglu}
\email{deniz.eroglu@khas.edu.tr}
\affiliation{Faculty of Engineering and Natural Sciences, Kadir Has University, 34083 Istanbul, Turkey}

\begin{abstract}
Reconstructing network dynamics from data is crucial for predicting the  changes in the dynamics of complex systems such as neuron networks; however, previous research has shown that the reconstruction is possible under strong constraints such as the need for lengthy data or small system size. Here, we present a recovery scheme blending theoretical model reduction and sparse recovery to identify the governing equations and the interactions of weakly coupled chaotic maps on complex networks, easing unrealistic constraints for real-world applications. Learning dynamics and connectivity lead to detecting critical transitions for parameter changes. We apply our technique to realistic neuronal systems with and without noise on a real mouse neocortex and artificial networks.
\end{abstract}

\maketitle

{\it Dynamical networks}, including power grids, food webs, climate networks, and neuron networks, described by dynamical units oscillating on complex networks, are fundamental components of our everyday lives. The ability to regulate network dynamics is crucial for predicting, thus, controlling these systems' behavior to acquire the desired functionality. Neuron networks are an important class of dynamical networks for human wellness since the changes in the interaction can lead to undesired pathological situations. For instance, epileptic seizures are associated with emergent neural network synchronization when the dynamical network parameters change \cite{Schindler2007}. Consequently, it is vital to anticipate critical transitions to neuronal synchronization and invent predictive technologies to detect early warning signals to prevent potential tragedies \cite{Eroglu2020}.
In the case of neuron network dynamics, consisting of intrinsic neuron function and the coupling scheme between neurons, the critical transitions to synchronization are not directly determinable. Therefore, the governing equation must be recovered from the observations of the nodes for forecasting the critical transitions due to parameter changes. 

The network dynamics reconstruction from data is a very active research field~\cite{Gao2022, wang2016data, nitzan2017revealing, stankovski2017coupling, Timme2014, casadiego2017model}. Various methods were proposed to infer the connectivity matrix under some constraints, such as the need for a system to be at steady-state~\cite{Gardner2003} or requiring prior knowledge about the dynamics~\cite{Timme2007, Shandilya2011, Yu2006} or the coupling strength~\cite{Ren2010}. In addition to the studies that reveal the connectivity matrix by control signals or analytical solutions, statistical learning approaches such as compressed sensing were also introduced to learn entire unknown dynamics~\cite{Wang2011, wang_2}, which also infers the connectivity structure. However, statistical learning techniques are not extendable for large networks or require long time series measurements. A natural question is then whether revealing the network dynamics of weakly interacting chaotic oscillators would be possible using relatively short data without requiring knowledge of the system's nodal behavior and coupling scheme. This question is especially relevant in weak coupling regimes, in which the synchronization regime is unstable and the decay of correlation is exponential for chaotic oscillators, meaning that similarity measures cannot capture the interaction topology.

This Letter reports a dynamical network reconstruction approach from time series observations by integrating mean-field approaches from dynamical systems theory with statistical learning tools. Neural networks are described by chaotic isolated dynamics \cite{Korn2003}, weakly interacting nodes \cite{Preyer2005} and interaction through scale-free type networks \cite{Werner2010}. Our reconstruction approach assumes that we have the mentioned neuroscientific setting and access to all nodes' data while the local dynamics of the nodes, the coupling function between them and the interaction structure are unknowns. Our methodology accurately identifies them using rather short time series and is independent from the network size, which is important since it is, generally, impossible to have long real-world observations, and real networks are large. Finally, as the reconstruction methodology includes mean-field approximations, the inferred model may not estimate the exact future states of the system due to the chaotic nature of the dynamical units. However, the reconstructed model allows us to predict the emergent collective behavior of dynamical networks considering parameter change, which is crucial to avoid undesired behaviors for real-world applications such as epilepsy seizures.

{\it Model.---} The network dynamics of weakly coupled and identical $n$ oscillators with interaction akin to diffusion is described by
\begin{equation}
\bm{x}_i(t+1)=\bm{f} (\bm{x}_i(t))+ \sum_{j=1}^{n} w_{ij} \bm{H}(\bm{x}_i(t),\bm{x}_j(t)) + \bm \eta_i(t)
\label{eq:netdynamics}
\end{equation} 
where $\bm x_i \in \mathbb{R}^m$, $\bm{f}\colon \mathbb{R}^m\to\mathbb{R}^m$ represents the isolated dynamics of nodes and we assume it is chaotic \cite{Barzel2013}. $\bm{H}$ is a diffusive coupling function ($\bm{H}(\bm x, \bm x)=\bm{H}(0)=0 \mbox{ and } \bm{H}(\bm x, \bm y) = -\bm{H}(\bm y, \bm x)$). $\bm W = [w_{ij}] \in \mathbb{R}^{n\times n}$ is the adjacency matrix of weighted and directed network where $w_{ij} \geq 0$ is the interaction strength from node-$j$ to node-$i$. The noise term, $\bm \eta_i(t)$, is uniformly distributed $\| \bm \eta_i(t) \| \le \eta_0$ for all nodes where $\eta_0$ is noise intensity. This network dynamics, Eq.~(\ref{eq:netdynamics}), is used to model numerous real-world applications including brain networks~\cite{izhikevich2007dynamical}, power grids~\cite{dorfler2013synchronization,motter2013spontaneous}, superconductors~\cite{watanabe1994constants}, and cardiac pacemaker cells~\cite{winfree2001geometry}.

{\it Reduction theorem.---} A low-dimensional reduction of Eq.~(\ref{eq:netdynamics}) is key for our network dynamics reconstruction approach. 
The reduction theorem applies a mean-field approach that relies on two main statements: (i) the statistical behavior of nodes' dynamics (frequency distribution of states) must be preserved and (ii) a large portion of nodes must be interacting with at least a few nodes in the network. These statements are satisfied with the given assumptions for the reduction theorem: \emph{chaotic local dynamics} of expanding maps and \emph{weak coupling} (to preserve the nodes' state distribution against fluctuations due to the interactions or external noise) and \emph{scale-free networks} (most of the nodes have small degrees $k\sim n^\epsilon$, and some nodes are hubs with degrees $k\sim n^{\frac{1}{2}+\epsilon}$ where $\epsilon$ is an arbitrarily small number), which also mimic brain network dynamics. Using the theorem, the coupling term of Eq.~(1) can be reduced as follows:

\begin{eqnarray}
\sum_{j=1}^{n} w_{ij} \bm{H}(\bm{x}_i,\bm{x}_j) &\approx& k_i \int \alpha \tilde{\bm{H}}(\bm x_i,\bm x_j) d  \mu (\bm x_j) \nonumber \\
&=& k_i \alpha (\tilde{\bm{V}}(\bm x_i) + \tilde C) = k_i \bm{V}(\bm x_i) + C\nonumber
\end{eqnarray}
where $\bm V$ is the effective coupling function, $k_i = \sum_j w_{ij}$ is the incoming degree of node $i$, $\bm{H} = \alpha \tilde{\bm{H}}$, $\alpha$ is a   multiplier for the coupling function, $\mu$ is a physical measure of the isolated dynamics, $C$ is the integration constant and the integral takes into account the cumulative effect of interactions on node-$i$. As the coupling term is reduced as a function of an invariant measure $\mu$, the reduction theorem works for a system in a steady state. Furthermore, to apply this mean-field approach-based reduction theorem, the statistical properties of individual dynamical systems must be preserved, which is satisfied by chaotic oscillators and weak coupling (see Supp.~Mat.~II and Ref.~\cite{Pereira2020}).
Then Eq.~(\ref{eq:netdynamics}) can be written as
\begin{eqnarray}
\bo x_i(t+1) = \bo f(\bo x_i(t)) + k_i \bm{V}(\bm x_i(t)) + C +  \bo \kappa_i(t) + \bo \eta_i(t)
\end{eqnarray}
where $\bm \kappa_i(t)$ is a small fluctuation for an interval of time that is exponentially large and depends on the state of neighbors of the $i$th node. 

\begin{figure}[!]
    \centering
    \includegraphics[width=.7\linewidth]{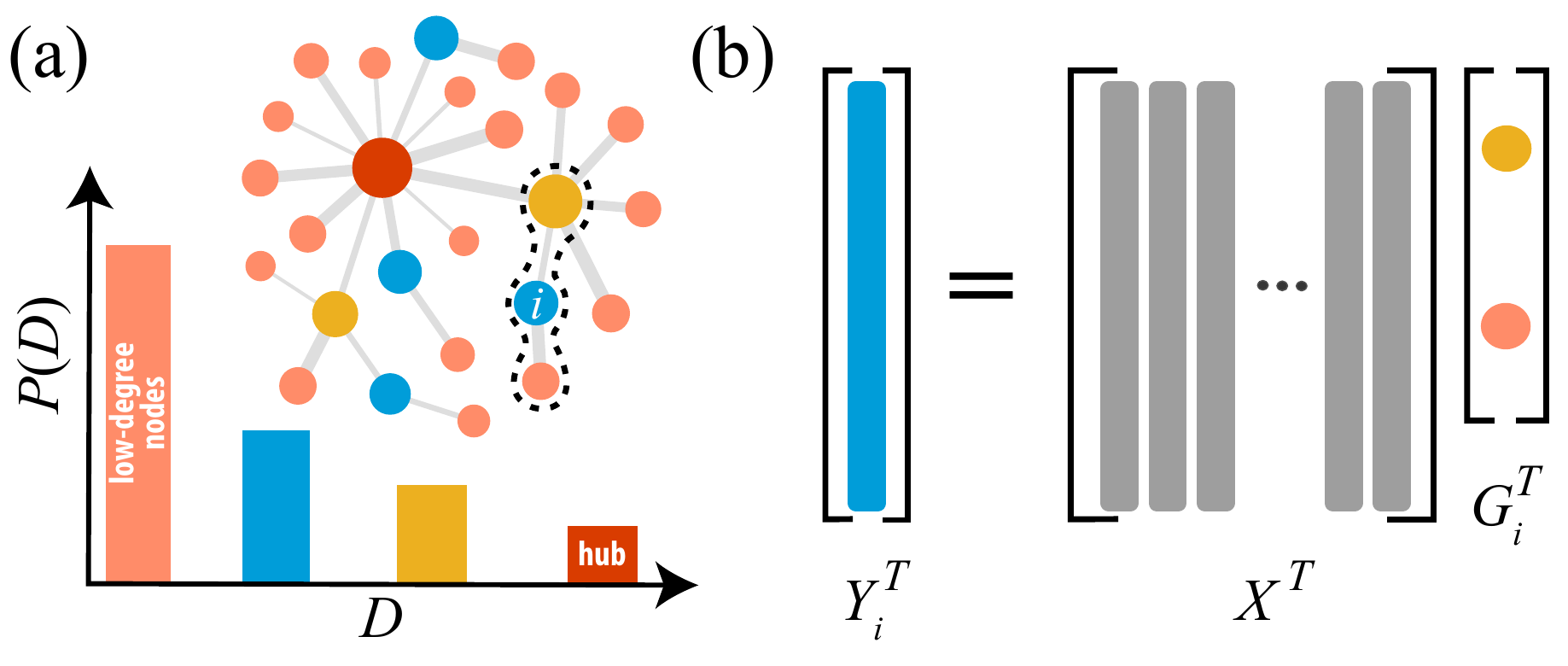}
    \caption{Illustration of the reconstruction scheme. (a) Performing sparse regression on each observation gives predicted models for each node. Nodes with the same in-degree are reconstructed with the same predicted models, which allow us to classify nodes concerning their in-degrees. As low-degree nodes are abundant and represent the isolated dynamics with a negligible noise, learning $\bm f$ is possible. Discarding the local dynamics, $\bm f$, from the hub's data gives the dominant coupling effect on the hub. Therefore, the coupling function $\bm H$ can be learned. (b) After learning $\bm f$ and $\bm H$, the problem is defined as a linear problem for each node by subtracting the local dynamics, which obtains the remaining interaction effect for each node. Sparse regression on the remaining interaction dynamics of node-$i$ where $i=1,\dots,n$ entirely reconstructs the dynamical networks. Nonzero $\bm{G}_i^T$ elements are the incoming connections for $i$th node (see Supp.~Mat.~Sec.~II for a step by step scheme).}
    \label{main_scheme}
\end{figure}

{\it Reconstruction scheme.---} To learn isolated dynamics $\bm f$ and coupling function $\bm H$, we first need to classify nodes regarding their degrees. According to the reduction theorem, nodes with a similar in-degree must have a similar governing equation. To identify the governing equations of each node independently from other nodes, we use sparse regression, particularly the Sparse Identification of Nonlinear Dynamical Systems (SINDy) technique \cite{Brunton2016}. We denote the data collection of node-$i$ by Eq.(\ref{eq:netdynamics}) as $\mathcal{X}_i = [\bm x_i(1), \dots,\bm x_i(T-1)]^T$ and $\mathcal{X}'_i = [\bm x_i(2), \dots,\bm x_i(T)]^T$. SINDy performs a sparse regression for the linear equation $\mathcal{X}'_i = \bm \Psi (\mathcal{X}_i) \bm \Xi_i$ to solve for $\bm \Xi_i = [\xi_i^1, \dots, \xi_i^p]^T$, which is a vector of coefficients that defines the dynamics, where $\bm \Psi = [\psi_1, \dots, \psi_p]$ represents a library of basis functions and is applied to $\bm x_i$ as
\[
\bm \Psi (\mathcal{X}_i)= \overset{\text{\normalsize candidate functions of $\bm x$}}{\left.\overrightarrow{\begin{bsmallmatrix}
    \psi_1 (\bm x_i(1)) & \psi_2 (\bm x_i(1)) & \dots & \psi_p(\bm x_i(1))\\
    \psi_1(\bm x_i(2)) & \psi_2(\bm x_i(2)) & \dots & \psi_p(\bm x_i(2))\\
    \vdots & \vdots & \ddots & \vdots \\
    \psi_1(\bm x_i(T-1)) & \psi_2(\bm x_i(T-1)) & \dots & \psi_p(\bm x_i(T-1)) 
\end{bsmallmatrix}}\right\downarrow}\begin{rotate}{270}\hspace{-.125in}time~~\end{rotate}
\]
where $p$ is the number of candidate functions in the library $\bm \Psi$. (For a detailed description for the basis library, see Supp.~Mat.~Sec.~VIII.) Sparse regression's goal is to determine the dynamics with a small number of functions in $\bm \Psi$ by finding active coefficients in $\bm \Xi_i$ (see Supp.~Mat.~Sec.~VII). Consequently, we obtain a predicted model for each node only using the associated node's own data, and we expect to learn similar models for the nodes with similar in-degree $k_i$. 
A distance matrix is obtained by normalized Euclidean distance to classify the predicted models $d_{ij}=(\Sigma_{k=1}^{p}~\frac{1}{V_k}~|\xi_i^k~-~\xi_j^k|^2 )^{1/2}$, where $|\cdot|$ is absolute value, $V_k$ is the variance of the predicted coefficients of the $k$th function in $\bm \Psi$. Assume $\bm \Xi_i$ and $\bm \Xi_j$ are two predicted models of nodes-$i$ and $j$, which are presented as a linear combination of some functions within the library $\bm \Psi$.  We get smaller $d_{ij}$ for a similar pair of nodes $i$ and $j$, while $d_{ij}$ will be large for distinct nodes, such as a low-degree node and a hub. An example computation of $d_{ij}$ can be found in Supp.~Mat.~Sec.~VIII. The histogram $P(D)$ is obtained by the row-sum of the distance matrix $D_i=\sum_j d_{ij}$, which provides an excellent classification for model similarities in terms of their degrees (Fig.~\ref{main_scheme}~(a)). The low-degree nodes are expected to be located in the highest bin of the histogram since the network has many low-degree nodes. The models recovered for low-degree nodes are determined as our $\bm f$ with a negligible fluctuation $\bm \kappa_i$. Contrarily, the distance $D_i$ is expected to be large for the hub nodes as they are the rarest. Therefore, hubs are located in the lowest bin of the histogram $P(D)$ (Fig.~\ref{main_scheme}~(a)). Note that the success of the reconstruction depends on the separability of the low-degree nodes and hubs with respect to their degrees, which means the network topology plays an important role here \cite{Eroglu2020}, which is further illustrated in Supp.~Mat.~Sec.~II.

We obtain the cumulative coupling effect on the hub by discarding learned isolated dynamics contribution from the identified hub node's data as $\mathcal{X}'_h-\bm f(\mathcal{X}_h)$ where $h$ denotes the hub node. We fit a function to the cumulative coupling effect on the hub and learn the coupling function $\bm H$ with a possible linear shift due to the integration constant $C$ (Supp.~Mat.~Sec.~II). The size of the linear shift can be easily estimated using $\bm H(\bm 0)= \bm 0$.  Inferring the interaction function is vital to reveal the network \cite{Liu2016}, and learning $\bm H$ from such reduced dynamics increases the feasibility of our approach. Introducing the Laplacian matrix, $\bm L$ with $L_{ij}=\delta_{ij}k_i-w_{ij}$ where $\delta_{ij}$ is the Kronecker delta ($\delta_{ii}=1$ and $\delta_{ij}=0$ if $i \ne j$) and assuming that $\bm H$ is a linear function, we can rewrite Eq.~(\ref{eq:netdynamics}) in a compact form as
\begin{equation}
\label{comp}
\bm X(t+1) = \bm{F}(\bm{X}(t)) - (\bm{L} \otimes \bm{H}) (\bm{X}(t)),
\end{equation}
where $\bm X=[\bm{x}_1,\cdots,\bm{x}_n]^T$, $\bm F(\bm X) = [\bm{f}(\bm{x}_1) , \cdots, \bm{f}( \bm{x}_n)]^T$ and $\otimes $ is the Kronecker product  (See Supp.~Mat.~Sec.~I for the derivation in terms of Laplacian matrix and \cite{pereira2013}). Defining $\bm Y(t) = \bm X(t+1) - \bm F(\bm X(t))$, Eq.~(\ref{comp}) can be written as $\bm{Y} = \bm{G} \bm X$ where $\bm G = -(\bm L \otimes \bm H)$. Finally, we complete the reconstruction by learning sparse matrix $\bm{G} \in \mathbb{R}^{mn \times mn}$, by solving the linear equation $\bm Y^T=\bm X^T \bm G^T$ using sparse regression, namely Least Absolute Shrinkage and Selection Operator (LASSO)\cite{Tibshirani1996}, as suggested in Ref.~\cite{Han2015}. LASSO adds the $\ell_1$ regularization penalty to the least-squares loss function to find the sparse coefficients (the links), and it is assessed as a compressed sensing approach \cite{candes2006stable}. Note that the linear equation can also be solved with $\ell_2$-norm for long time series, however, as we are interested in short data (in the case of the length of the time series  $T < mn$) the compressed sensing approach must be employed \cite{candes2006stable}. Consequently, we learn the connectivity matrices $\bm G$ and $\bm L$ as seen in Fig.~\ref{main_scheme}~(b). It is also important to note that learning the equations of all nodes by a single sparse regression without the reduction theorem is only possible for relatively small networks. The library extension, due to network size, causes a statistically correlated data matrix that quickly fails on reconstruction \cite{Novaes2021, desilva2020pysindy}. A detailed discussion can be found in Supp.~Mat.~Sec.~XII.

{\it Mouse neocortex reconstruction.---} 
A weighted and directed neural network ($987$ nodes and $1536$ edges), representing a mouse neocortex~\cite{Kasthuri2015, Vogelstein2018}, is considered (Supp.~Mat.~Sec.~IV). To mimic neurons, we used electrically coupled Rulkov maps akin to diffusion as
\begin{eqnarray}
u_i(t+1) &=& \frac{\beta}{1+u_i(t)^2} + v_i(t) - \sum_j L_{ij} u_j + \eta_i \nonumber \\ 
v_i(t+1) &=& v_i(t) - \nu u_i(t) - \sigma \nonumber
\end{eqnarray}
where the fast variable $u_i$ is the membrane potential and the slow variable $v_i$ is the ion concentration variation \cite{rulkov-variables}. The constant parameters $\beta=4.1$ and $\nu=\sigma=0.001$ are fixed for chaotic bursting dynamics~\cite{Rulkov2002}.
\begin{figure}[!]
    \centering
    \includegraphics[width=.7\linewidth]{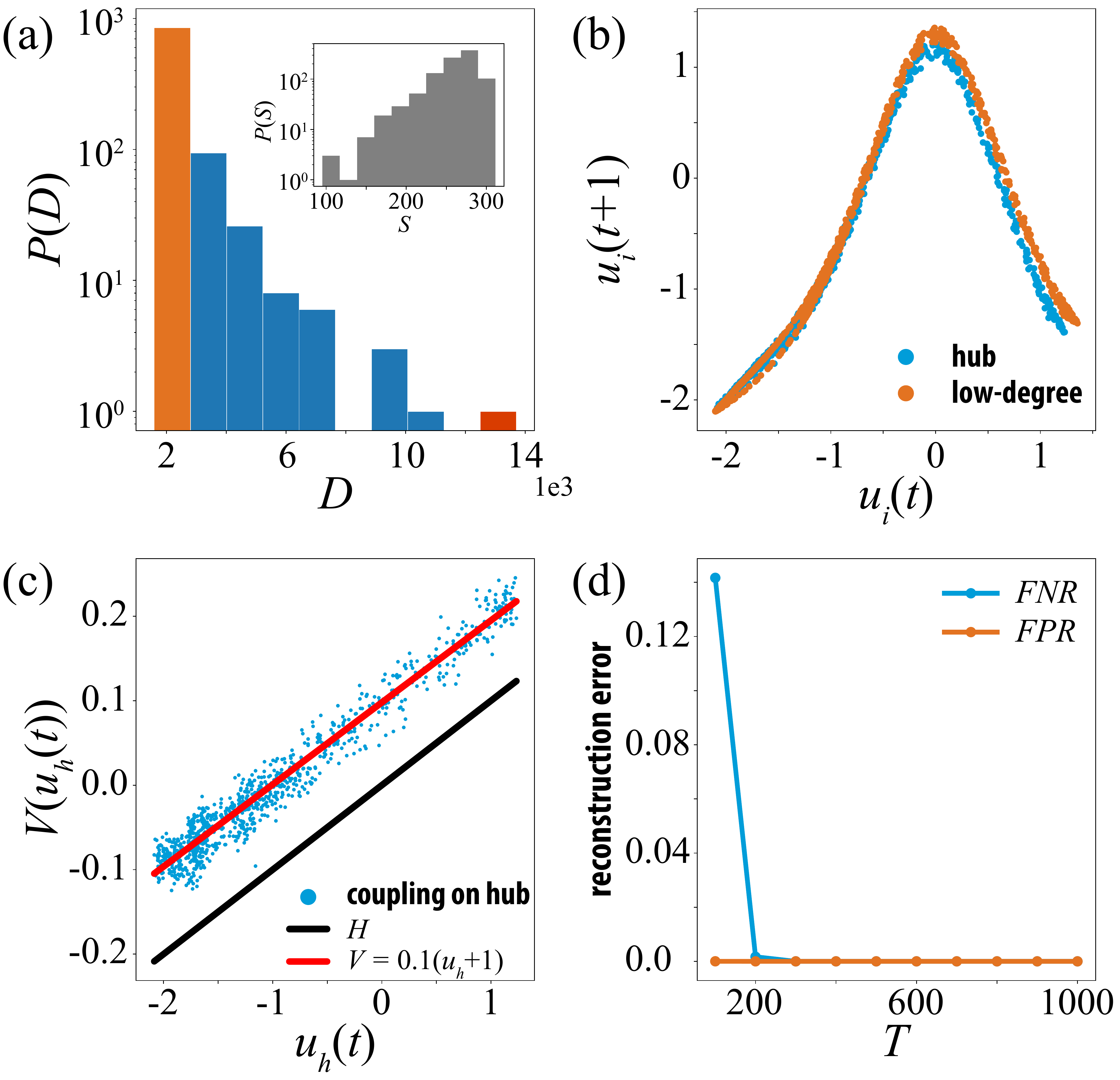}
    \caption{Reconstruction procedure for weakly electrically coupled Rulkov maps on a real mouse neocortex network. (a) Node-similarity histogram determines the low-degree nodes (orange bar) and the hub (red bar). Inset: Histogram $P(S)$ presents the correlations between the original time series. It is impossible to infer the connectivity structure from the correlations due to the chaotic nature of the Rulkov maps. (b) The return maps of a low-degree node and the hub are slightly different due to weak coupling effect. (c) The effective coupling, $V(u)$, shifts through the horizontal direction due to the integral constant $1$. (d) $FNR$ for different lengths of time series. $FPR$ are zero for all time series lengths.}
    \label{rulkov_results}
\end{figure}

\begin{figure}[!]
    \centering
    \includegraphics[width=.7\linewidth]{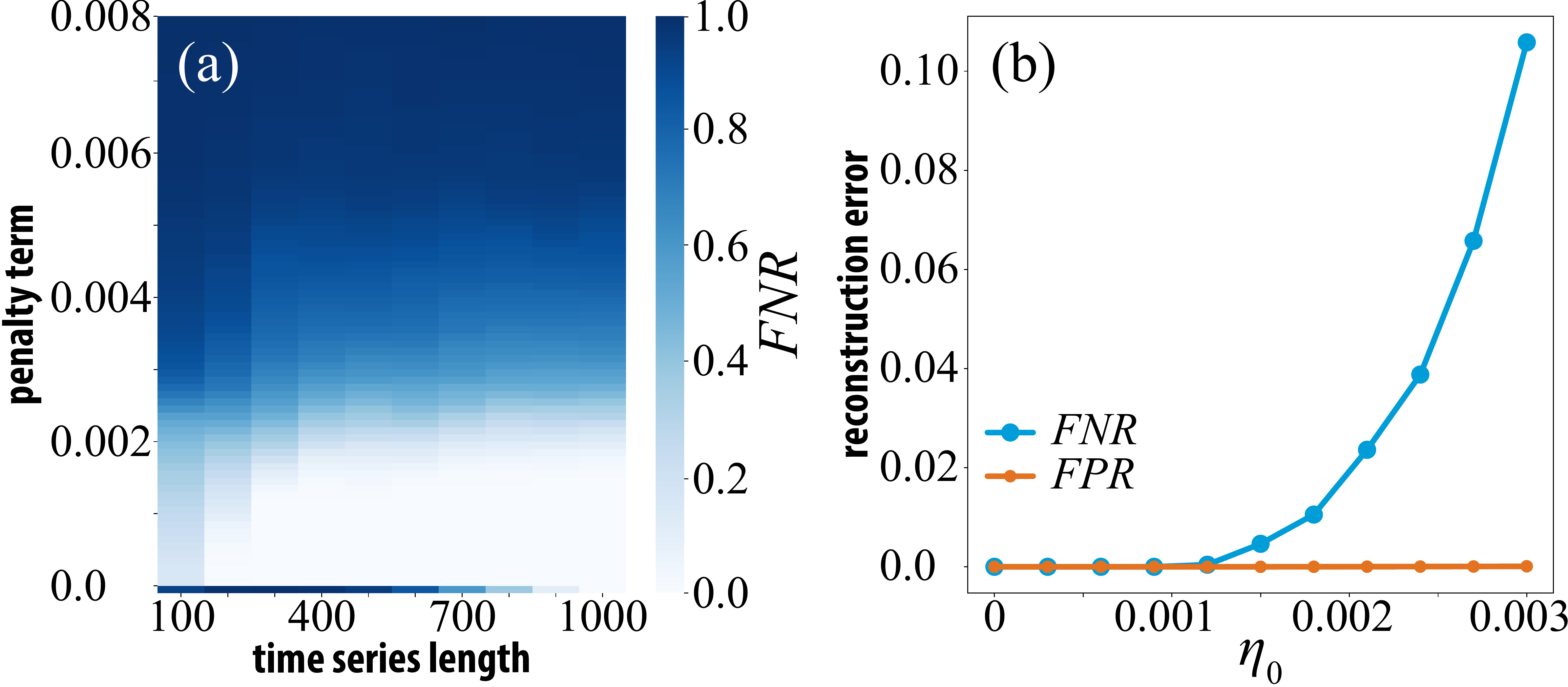}
    \caption{(a) The reconstruction performance for different time series lengths and a series of penalty terms as $FNR$. $FPR$ always equals 0 for this case. (b) Noise effect on reconstruction performance on real network as $FNR$ and $FPR$.}
    \label{real-network}
\end{figure}

{\it Noise-free case.---}
Following the reconstruction scheme, the nodes are classified using the similarity histogram Fig.~\ref{rulkov_results}(a) for $\eta_0 = 0$, while the pairwise Pearson correlations do not show any information about the degrees (inset in Fig.~\ref{rulkov_results}(a)).
The difference between return maps of a low-degree node and the hub is illustrated in Fig.~\ref{rulkov_results}(b), and the comparison against an isolated node dynamics is given in Supp.~Mat.~Sec.~V-A.
Effective coupling ${\bm V}(\bm x)$ is found approximately $[0.1(u + 1), 0]$, meaning that $\alpha$ is $0.1$ and the linear shift is $1$ on $u$-variable due to $\bm H( \bm 0) = \bm 0$ (Fig.~\ref{rulkov_results}(c)).  Finally, we learn the network topology by solving the linear equation $\bm{Y} = \bm{G} \bm X$ using the learned $\bm f$ and $\tilde{\bm H}$.
We measure the reconstruction error using the fraction of the false negatives (positives) out of the positives (negatives), $FNR$ ($FPR$). The $FNR$ ($FPR$) equals 0 for perfect reconstruction (Supp.~Mat.~Sec.~III). 
Here, we use the ground truth Laplacian matrix to assess the accuracy of the reconstruction.
When the ground truth is not available, the learned model can be evaluated by the cross-validation techniques (Supp.~Mat.~Sec.~XI). 
The reconstruction error is found to be almost zero for a data length larger than $T>200$ (Fig.~\ref{rulkov_results}(d)). 
Furthermore, a systematical evaluation of our approach according to penalty terms and time series lengths is performed. When the number of nodes $n$ times the dimension of the local dynamics $m$ exceeds the time series length $nm > T$, it corresponds to an underdetermined linear problem. Even for short data ($T \approx 200$),  a successful reconstruction is possible  with a small penalty term when $mn = 1974$ (Fig.~\ref{real-network}(a)). Note that if the problem is overdetermined $nm < T$, then $\ell_2$-norm regression can be used for faster computations (Supp.~Mat.~Sec.~IX).
Furthermore, we provide reconstruction analysis for H\'enon map and the Tinkerbell map in Supp.~Mat.~Sec.~V. We also perform our procedure on a macaque monkey visual cortex network \cite{rhesus,Vogelstein2018}. This network is not scale-free; however, assuming the local dynamics and the hub node are known, we reconstructed the network dynamics (see Supp.~Mat.~Sec.~VI).

\begin{figure}[!]
    \centering
    \includegraphics[width=0.7\linewidth]{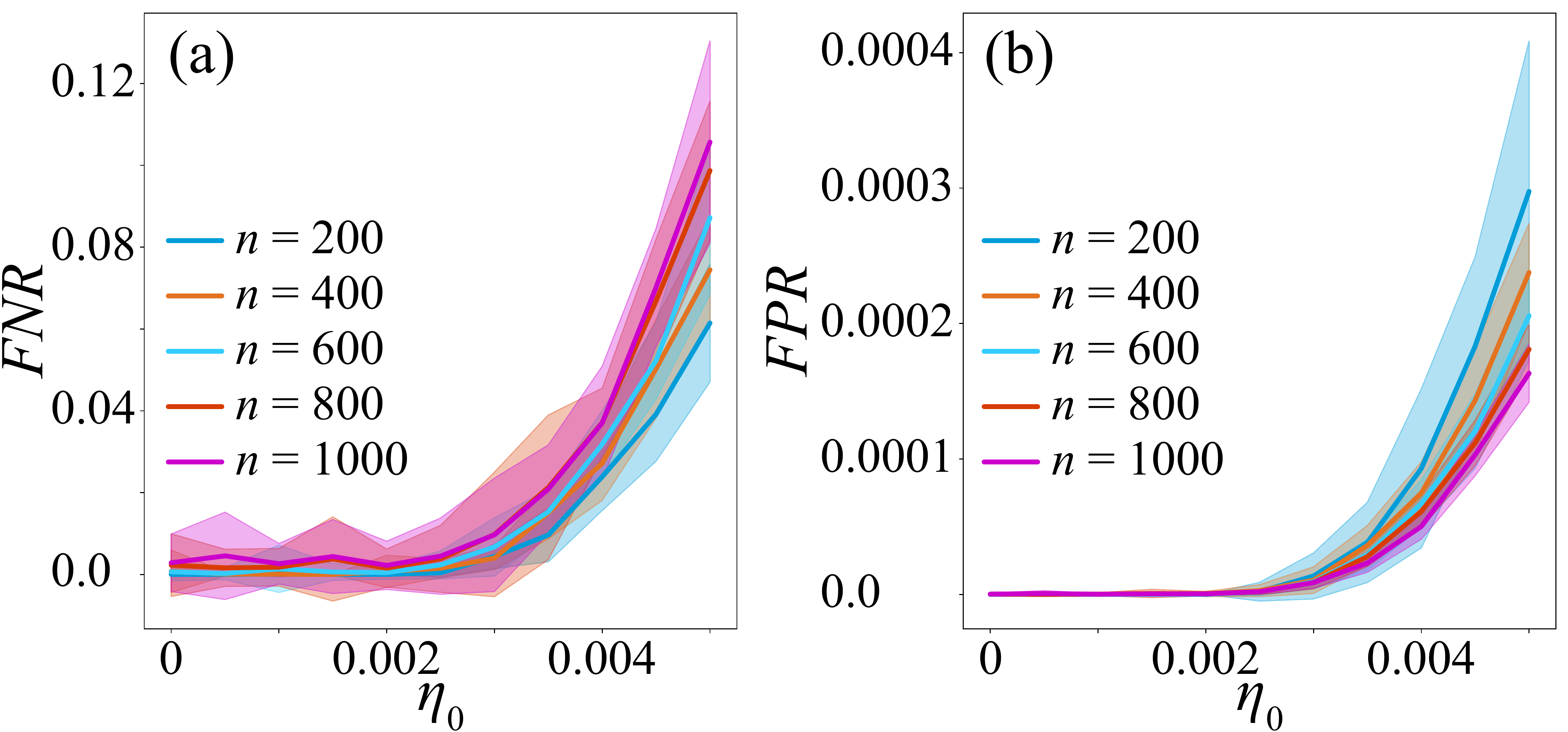}
    \caption{Average noise effect on reconstruction performance illustrated using (a) $FNR$ and (b) $FPR$ for $50$ realizations of simulations using random scale-free networks of system sizes $n=200, 400, 600$, $800$ and $1000$. Time series length is fixed as $500$ during all the simulations. The shaded regions represent the corresponding standard deviations.}
    \label{synthetic-networks}
\end{figure}

{\it Noise effect on reconstruction performance.---} To measure the robustness of our methodology against noise, we systematically perform the reconstruction procedure for various $\eta_0$ values (Eq.~(\ref{eq:netdynamics})) on synthetic as well as real mouse neocortex networks. The reconstruction approach is robust to small noise intensities for the real-world example (Fig.~\ref{real-network}(b)). The average robustness of the reconstruction procedure over $50$ different directed and weighted random scale-free networks is given in Fig.~\ref{synthetic-networks}. 
The scale-free networks are generated using the algorithm in Ref.~\cite{directed_scale-free}, and weights are assigned uniformly from the interval $[0.8,1.2]$. The algorithm generates undesired self-loops and multiple edges, so first we remove them. As the system size grows, the reconstruction performance reduces with respect to $FNR$ for increasing noise intensity $\eta_0$ (Fig.~\ref{synthetic-networks}(a)), since the noise becomes more dominant than the weak coupling, which prohibits learning $\bm H$ using the coupling effect. As similar to real-world application, $FPR$ results are also negligibly small for the noise induced reconstruction case (Fig.~\ref{synthetic-networks}(b)). We also performed experiments by generating denser networks where the reconstruction technique fails, therefore, the network sparsity is crucial for the reconstruction (see Supp.~Mat.~Sec.~X).

{\it Prediction of emergent behavior.---} Although we recover the true network structure, the reduction theory approximates the isolated dynamics $\bm f$ and coupling function $\bm H$ with small bounded fluctuations. As the local dynamics is chaotic and we consider possible noise effect, time evolution forecasting of the given system can diverge from the original system. However, as we recovered the network dynamics with high accuracy for the given data as, 
\begin{equation}
\bm{X}(t+1) = \bm{F}(\bm{X}(t)) + \gamma \bm{G}\bm{X}(t),
\end{equation}
where $\gamma$ is the coupling control parameter and initially it is $\gamma=1$ for the reconstructed model. Then, it is possible to detect critical coupling strength factor $\gamma_c$ to predict the emergent behavior of the dynamical network by fully analytical techniques if the reconstructed coupling function is identity matrix \cite{eroglu2017synchronisation}. For general coupling functions, master stability function \cite{Pecora_msf} or connection graph method \cite{Belykh2004} can be performed for the detection (see Supp.~Mat.~Sec.~XIII).

{\it Conclusions.---} The network is fully recovered by our approach in the setting where $\bm f$ is chaotic, the network is scale-free and the coupling is weak. Because of the weak coupling between nodes and the chaotic nature of the local dynamics, the correlation between measured time series decays exponentially; therefore, it is impossible to reconstruct such complex systems by conventional methods. Cutting-edge autonomous statistical learning techniques also fail when the network size is large. The key idea in our procedure is splitting the model equation into parts by reduction theorem and inferring each unknown ($\bm f, \bm H$ and $\bm W$) one by one using sparse recovery. Although the reduction theorem is not established for the general chaotic discrete maps, we showed its validity on various maps. Our  approach guarantees the full reconstruction for the noise-free case and small noise intensity even for relatively short time series with no limitations on the network size.
However, the quality of reconstruction decreases for increasing noise and the destructive effect of the noise also increases with increasing system size. Finally, obtaining the network dynamics allows one to predict the emergent behavior under parameter changes. There are  available regression-based approaches that can learn network topology using short time series~\cite{casadiego2017model}; however, it is impossible to detect the critical transitions with only the connectivity. The ability to detect such transitions is crucial for applications such as a transition to collective behavior in the brain network, which can lead to undesired implications. Thus, it is desirable to put forward precautionary norms to avert potential disasters.

{\it Data and code availability.---}The data we used in this study can be regenerated by running the code which is publicly available on GitHub~\cite{github}.

We are indebted to Tiago Pereira, Matteo Tanzi, Sajjad Bakrani, Arash Rezaeinazhad, Thomas Peron and Jeroen Lamb for enlightening discussions. This work is supported by The Scientific and Technological Research
Council of Turkey (TUBITAK) under Grant No. 118C236. D.E. acknowledges support from the BAGEP Award of the Science Academy.

\let\oldaddcontentsline\addcontentsline
\renewcommand{\addcontentsline}[3]{}

\let\addcontentsline\oldaddcontentsline
\end{bibunit}

\newpage

\onecolumngrid

\setcounter{page}{1}

\begin{center}
{\large \bf Reconstructing network dynamics of coupled discrete chaotic units from data:} Supplemental Material  \\
\vspace{0.3cm}
{ Irem Topal and Deniz Eroglu}\\
\vspace{0.1cm}
{\small
Faculty of Engineering and Natural Sciences, Kadir Has University, 34083 Istanbul, Turkey}
\end{center}

\tableofcontents

\newpage

\section{Connectivity matrix in terms of Laplacian}
Network dynamics with diffusive interaction is given by,
\begin{equation}
\bm{x}_i(t+1)=\bm{f} (\bm{x}_i(t))+ \sum_{j=1}^{n} w_{ij} [\bm{H}(\bm{x}_j) - \bm{H}(\bm{x}_i)]
\label{eq:netdynamics}
\end{equation} 
where $\bm x \in \mathbb{R}^m$, $\bm{f}\colon \mathbb{R}^m\to\mathbb{R}^m$ is chaotic dynamics of isolated nodes and $\bm H$ is a diffusive coupling function \cite{Barzel2013}. $\bm W = [w_{ij}] \in \mathbb{R}^{n \times n}$ is a weighted and directed adjacency matrix. Diffusive nature of the interaction allows represent the coupling in terms of the Laplacian matrix.
\begin{eqnarray}
\sum_{j=1}^{n} w_{ij}[\bm{H}(\bm{x}_j) - \bm{H}(\bm{x}_i)] &=& \sum_{j=1}^{n} w_{ij}\bm{H}(\bm{x}_j) - \bm{H}(\bm{x}_i) \sum_{j=1}^{n} w_{ij} \\ \nonumber
&=& \sum_{j=1}^{n} w_{ij}\bm{H}(\bm{x}_j) - k_i \bm{H}(\bm{x}_i)\\ \nonumber
&=&\sum_{j=1}^{n} (w_{ij} - \delta_{ij} k_i )\bm{H}(\bm{x}_j) \end{eqnarray}
where $k_i = \sum_j w_{ij}$ is the incoming degree of node $i$ and $\delta_{ij}$ is the Kronecker delta. By introducing the Laplacian matrix, $\bm L$ with $L_{ij}=\delta_{ij}k_i-w_{ij}$ we obtain,
\begin{equation}
\bm{x}_i(t+1)=\bm{f} (\bm{x}_i(t)) - \sum_{j=1}^{n} L_{ij}\bm{H}(\bm{x}_j)
\label{laplacian-form}
\end{equation}
Then we can rewrite Eq.~(\ref{laplacian-form}) in a compact form as,
\begin{equation}
\label{comp}
\bm X(t+1) = \bm{F}(\bm{X}(t)) - (\bm{L} \otimes \bm{H}) (\bm{X}(t)),
\end{equation}
where $\bm X=[\bm{x}_1,\cdots,\bm{x}_n]^T$, $\bm F(\bm X) = [\bm{f}(\bm{x}_1) , \cdots, \bm{f}( \bm{x}_n)]^T$  and $\otimes $ is the Kronecker product \cite{pereira2013}. 

\section{Network Reconstruction Scheme on a toy model} 
\begin{figure}
\centering
\includegraphics[width=0.8\linewidth]{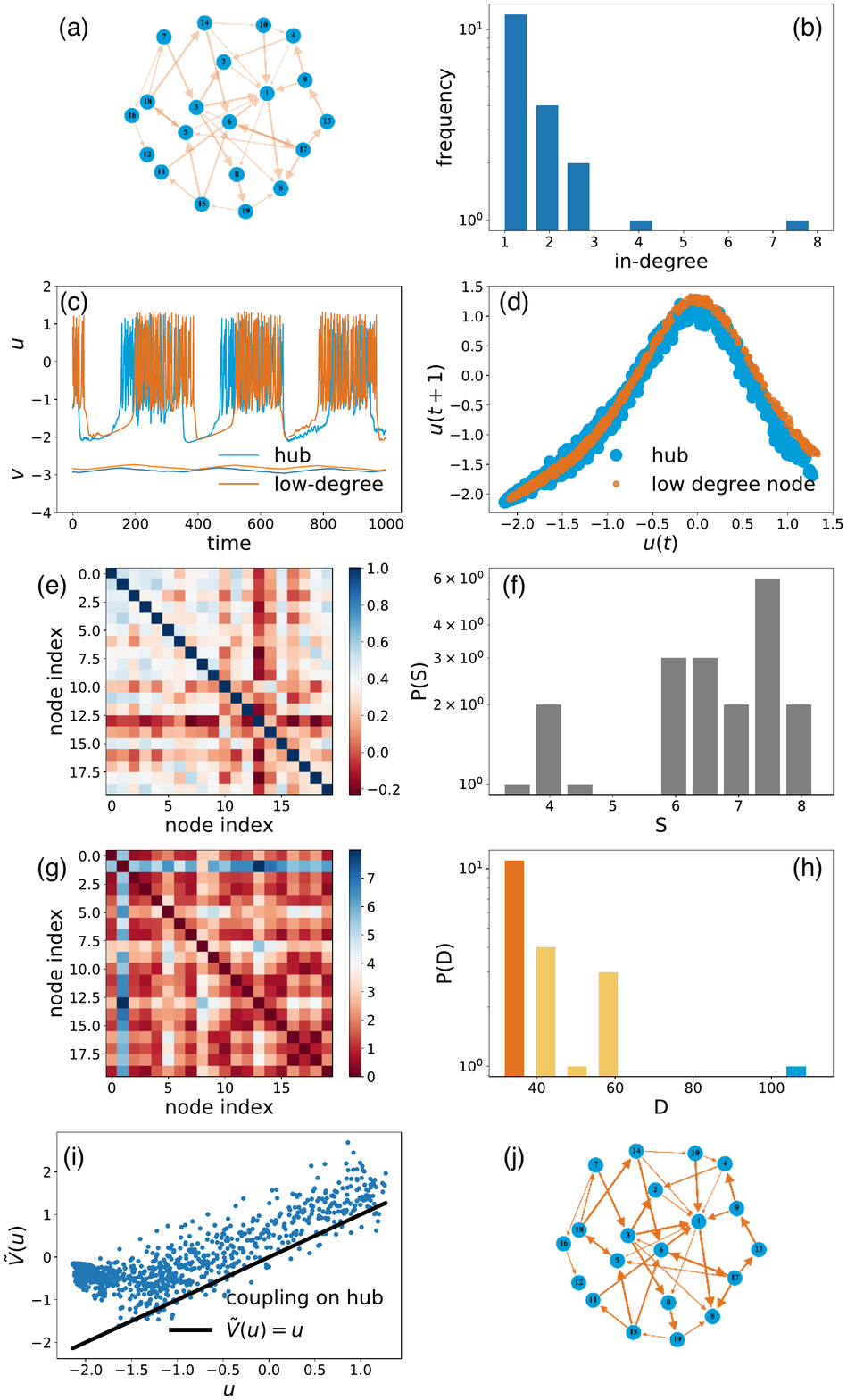}
    \caption{(a) A directed and weighted network. As the network is unknown initially, its links were illustrated indistinctly. (b) Corresponding in-degree distribution of the network. (c) Time series of a low-degree node and the hub. (d) The return maps of these time series, the hub's data is dispersed according to the low-degree node due to the coupling effect accumulated on the hub. (e) Pairwise Pearson correlation matrix of the time series, which shows no similarity between, even connected, nodes dynamics due to the chaos. (f) The histogram obtained by the row-sum of the time series similarity matrix (subplot (e)) does not reflect any information about the network's original degree distribution; therefore, it is impossible to identify a low-degree node or the hub. (g) Pairwise distances between predicted models were obtained in the first step of our reconstruction procedure. Nodes with similar in-degrees show high similarity. (h) The histogram of the row-sum of the distance matrix (subplot (g)), which does also not contain the exact in-degree distribution of the network (subplot (b)), however, allows for the low-degree and hub node identification. The abundant models (the highest bin in (h)) represent the low-degree nodes, and these models represent the local dynamics, $\bm f$.  The most distinct model (the lowest bin in (h)) belongs to the hub node. We learn the coupling function by filtering the local dynamics from the hub data. (i) The remaining coupling effect on the hub is plotted, and the reduction theorem recovers the coupling function. (j) After recovering the local dynamics and the coupling function, we solve the linear equation $y=Ax$ by sparse regression to obtain the network structure.}
    \label{scheme}
\end{figure}

We present a step-by-step reconstruction algorithm using weakly coupled chaotic maps on a directed and weighted network of size $n = 20$ Fig.~\ref{scheme}(a). The network is heterogeneous; the low-degree nodes are abundant, and the hub is the rarest, as illustrated in the network's in-degree distribution Fig.~\ref{scheme}(b). The network has Rulkov maps interacting with each other through their $u$-components diffusively: 
\begin{eqnarray}
u_i(t+1) &=& \frac{\beta}{1+u_i(t)^2} + v_i(t) - \sum_{j=1}^{n} L_{ij} u_j(t)\\ 
v_i(t+1) &=& u_i(t) - \nu u_i(t) - \sigma \nonumber
\label{governing_equations}
\end{eqnarray}
where the constant parameters $\beta=4.1$ and $\nu=\sigma=0.001$ are fixed for chaotic bursting dynamics. We first simulate the system $15000$ time steps and discard the first $14000$ steps as a transient. We should note that the local dynamics (Rulkov map), the coupling function ($u$-coupling) and the connectivity matrix ($L_{ij}$) are all unknowns during the procedure; however, we present every piece of information in this toy model to show the method rigorously. 

Fig.~\ref{scheme}(c) shows the $2$-dimensional data coming from the hub and one of the low degree nodes. We present the return maps of the same nodes in Fig.~\ref{scheme}(d), the dispersion of the hub node due to the dominant coupling effect can be seen. On the other hand, we compute the pairwise Pearson correlation coefficients $s_{ij}$ between the time series, see Fig.~\ref{scheme}(e). Then we obtain the histogram $P(S)$, Fig.~\ref{scheme}(f), by the row-sum of the correlation matrix in Fig.~\ref{scheme}(e). Given the chaotic nature of the Rulkov maps, correlations between time series does not give any information about the network connectivity.

The reconstruction procedure is given by the following steps:
\begin{enumerate}
\item{\bf Classification of the nodes.~}

The separability of low-degree nodes and hubs is crucial since the reconstruction procedure is not applicable if we cannot classify which nodes are low-degree or hubs at the beginning of the reconstruction recipe. The classification leads us to learn the isolated dynamics $\bm f$ and the coupling function $\bm H$ as follows: we, first, assume that the coupling is \emph{weak}, meaning that if we classify the low-degree nodes, we can consider their dynamics as local dynamics $\bm f$ with some small fluctuations. If we classify which node is the hub and $\bm f$, then we can learn the coupling function $\bm H$ by discarding $\bm f(\bm x_h(t))$ from the hub's time series $\bm x_h(t+1)$ using the cumulative effect of interactions. Therefore, our first aim is to classify the nodes regarding their degrees. 

For the classification, we first learn a governing equation for each time series using sparse regression methods (Supp.~Mat.~Sec.~\ref{sec:regression}) with a basis composed of candidate functions (Supp.~Mat.~Sec.~\ref{sec:library}). 
After obtaining the $n$ predicted models, we measure the difference between the models using normalized Euclidean distance between the coefficients of predicted models, $d_{ij}$. A distance matrix $D = [d_{ij}]$ summarizes the mismatches between the predicted models.  Fig.~\ref{scheme}(g). The row-sum of the distance matrix $D_i=\sum_j d_{ij}$ gives a histogram $P(D)$, and each bin of the histogram contains nodes with a similar degree (Fig.~\ref{scheme}(h)). 

It is also important to note that weak coupling is also crucial to avoid possible synchronization in the network dynamics. If the system is synchronized, the time series will be identical; in other words, all the predicted models will be identical. This will make the classification impossible. Therefore, the weak coupling is quite crucial for this algorithm to work as desired.

\item{\bf Learning the local dynamics.~}
To mimic brain network dynamics, we assume that the connectivity matrix $\bm W = [w_{ij}]$ represents a scale-free network whose degree distribution follows a power law; in other words, most of the nodes have small degrees $k\sim n^\epsilon$, and some nodes are hubs with degrees $k\sim n^{\frac{1}{2}+\epsilon}$ where $\epsilon$ is an arbitrarily small number. Therefore,  the highest bin of the histogram (the most abundant frequencies) contains the low-degree nodes' models since there are many low-degree nodes (the orange bar in Fig.~\ref{scheme}~(h)). As it is likely to have nodes having incoming links, $k_i > 0$, in the set of low-degree nodes, the models can contain some small fluctuations. Thus, we average the predicted coefficients of the low-degree nodes' models to infer the local dynamics as accurately as possible. 

\item{\bf Learning the coupling function.~}
Using the scale-free network topology assumption as in the previous section, we understand that the lowest bin of the histogram (the rarest frequencies ) contains the most different models, which belong to the hubs (the blue bar in Fig.~\ref{scheme}(h)). Then discarding the local dynamics of the hub ($\bm f( \bm x_h(t)$) from the hub measurement ($\bm x_h(t+1)$) gives the dominant coupling effect, which we aim to learn. Fig.~\ref{scheme}(i) presents the remaining coupling effect. We fit a function to the data and learn an approximate coupling function plus an integration constant arising from the reduction theorem. To be more precise with the integration constant, for instance, if the diffusive coupling function is given as $h(x,y) = \phi (y) - \phi (x)$
then the effective coupling $V$ for the expending maps is $V = \int h(x,y) d\mu(y) = -\phi (x) + C$
where $C$ is the integration constant, which we call a \emph{possible linear shift} for the numerically recovered effective coupling $V$. This is easy to estimate from the fitted coupling (Fig.~\ref{scheme}~(i)) since we assume that $\phi(0)=0$. Therefore, the appearing linear shift is the integration constant $C$, and its value can vary for different dynamical networks. Furthermore, the effective coupling is an approximation, which becomes exact when the network model is considered as a directed star graph with the link directions from low-degree nodes to the hub at the thermodynamical limits as $n \to \infty$. We have an approximated effective coupling for tree-like graphs, as scale-free graphs, under our assumptions related to node degrees (mentioned above). The size of the bounded fluctuation term $\kappa_i$ in the reduction (Eq.~(2) in the main text) can be considered as the measure of approximation. 

Here, the weak coupling and chaotic dynamics assumptions take an important role again for more general networks to preserve the statistical behavior of the nodes' dynamics. Although using chaotic oscillators helps to keep the state distribution of the states stable, this property can be destroyed by a large coupling. Therefore, if the network is denser, then the coupling constant should be scaled to a smaller coupling strength for our reduction theorem to work. For instance, if the hubs are rare in a scale-free network, the coupling strength must be scaled by $n^{1/2}$. If all the nodes are hubs, in other words, if the network is all-to-all connected, then the coupling strength must be scaled by $n$.  In our work, we scaled the weights in $\bm W$ with $n^{1/2}$, as we considered scale-free networks (containing many low-degree nodes and rare hubs) since the work focuses on reconstructing brain networks. It is also important to remind that, for the node classification step, the synchronization of the nodes disallows the reconstruction since the machine learning-based methods cannot distinguish the data source between synchronized oscillators. Consequently, the chaotic oscillators also play an important role here since if we have periodic and identical oscillator networks, any positive coupling strength can synchronize the model, which is undesired for the reconstruction part. Therefore, the weak coupling strength and the chaotic dynamical regime are crucial assumptions for the theory and also for the numerical steps to learn the coupling function $\bm H$. 
 
\item{\bf Learning the connectivity matrix.~} 
After learning the functions $\bm f$ and $\bm H$, we define $\bm Y = \bm X(t+1) - \bm F(\bm X(t))$, so Eq.~(\ref{comp}) can be written as $\bm{Y} = \bm{G} \bm X$ where $\bm G = -(\bm L \otimes \bm H)$. Therefore, learning the links becomes a  regression problem by solving the linear equation $\bm Y^T=\bm X^T \bm G^T$. Here we employ a compressed sensing approach using $\ell_1$-norm, called LASSO, to solve the equation since we are interested in reconstructing large network dynamics from short data. In other words, we aim to solve the problem when the length of the time series  $< m \times n$ where $m$ is the dimension of the local dynamics and $n$ is the system size. The problems under this condition are called underdetermined regression problems. Finally, we identify the exact Laplacian matrix of the network by solving the linear equation (Fig.~1 (j)). 

Note that, in the previous steps of the procedure, the potential interaction functions of two nodes (or higher-order interactions if the hypernetworks are interested) are not involved in the basis library. Meaning that the size of the library is not grown exponentially. Therefore, the algorithm does not require a longer time series to reconstruct the network dynamics as the network size, $n$, is increasing. In this last step, we also only sparsely solve a linear equation where $\bm G \in \mathbb{R}^{mn \times mn}$, which is a square matrix. The discussion on the relation between the length of the time series and the sparsity can be found in Supp.~Mat.~Sec.~\ref{sec:whyreduction}.

 \end{enumerate}
 
\section{Assessing the performance}
We use two standard indices to assess the reconstruction performance of our procedure: the fraction of the false negatives out of the positives ($FNR$) and the fraction of false positives out of the negatives ($FPR$). These metrics correspond to two types of errors in network reconstruction: missing a link where a link indeed exists (false negatives) or assigning a link between two nodes when they are not connected (false positives): 
\begin{equation}
 FNR = \frac{FN}{P} \text{ and }FPR = \frac{FP}{N} \nonumber
\end{equation}
where $FN$ evaluates errors on the non-zero terms ($P$) of the correct matrix (underestimation) and $FP$ evaluates errors of the zero terms ($N$) (overestimation). We count $FN$ and $FP$ links by $\sum_{i,j}\Theta(|L_{ij} - \hat L_{ij}| > \epsilon)$ to assess the accuracy of the link strengths where $\hat L_{ij}$ is the predicted Laplacian and $\epsilon$ is a tolerance value set as $10$ times smaller than the smallest nonzero value in $\bm L$. We use a fixed tolerance value as $\epsilon = 0.0001$ for  all analyses.
\newpage
\section{Mouse neo-cortex network}
\label{real-net}
\begin{figure}[H]
    \centering
    \includegraphics[width=.8\linewidth]{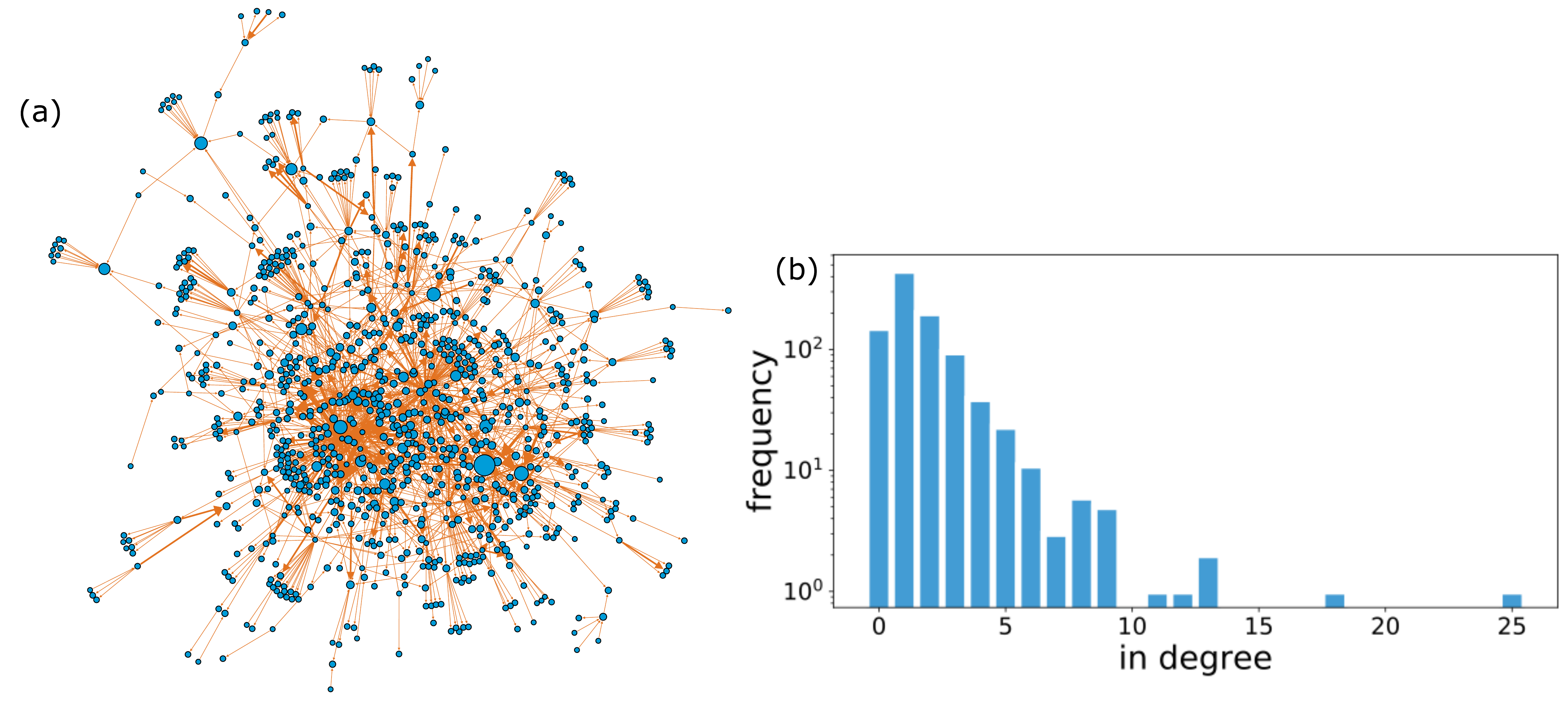}
    \caption{(a) Pre-processed directed and weighted mouse neo-cortex network. The sizes of the nodes denote their in-coming degrees and the thickness of the arrows denote their weights. (b) In-degree distribution of the associated real network.}
\label{kasthuri}
\end{figure}

We choose a neuronal connectivity consisting $1029$ nodes representing a small volume of a young adult mouse neocortex \cite{Kasthuri2015}. The network is directed and contains multi-edges between some nodes. We considered the number of these multi edges between nodes as weights. Furthermore, the network contains some disconnected parts, we removed them and kept only the giant component. The resulting weighted and directed network has $n=987$ nodes and $m=1536$ links, which is illustrated in Fig.~\ref{kasthuri}~(a) with its in-degree distribution in Fig.~\ref{kasthuri}~(b). We downscale entries of the connectivity matrix by a small constant $0.1$, and then normalize by in-degree of hub, $k_h$, so that the network is not synchronized. 

\section{Applications on mouse neo-cortex network}

This section presents the results of the reconstruction procedure on various dynamical systems, namely Rulkov map, H\'enon map and Tinkerbell map. For all the systems, we use mouse neo-cortex network which described in Section~\ref{real-net}.

\subsection{Rulkov map}
We have already presented Rulkov map results in the main text. In our reconstruction scheme, we assume that a low-degree node can be taken as an isolated node to estimate local dynamics $\bm f$. Here, we compare the return maps of a low-degree node and a hub node's dynamics with an isolated Rulkov map under the effect of weak coupling (Fig~\ref{return_maps_rulkov}).

\begin{minipage}[H]{\linewidth}
    \centering
    \includegraphics[width=0.6\linewidth]{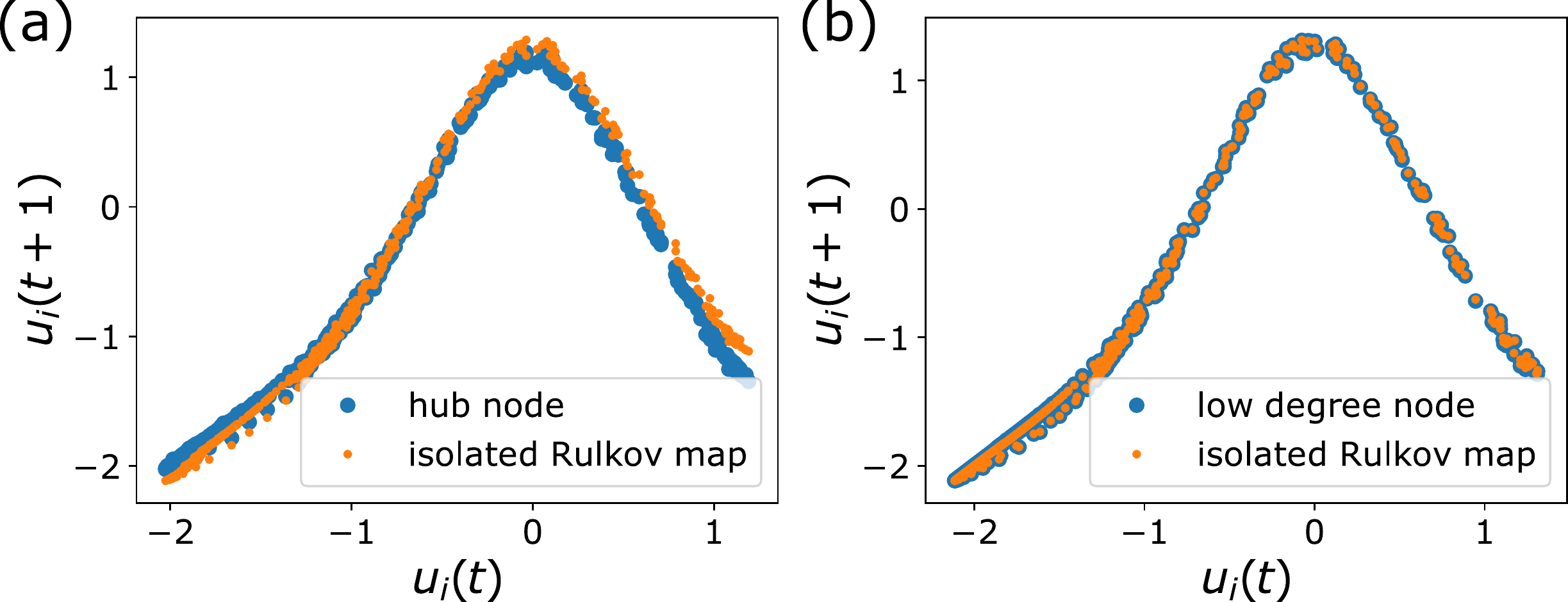}
    \captionof{figure}{a) The return map of the hub's data and isolated Rulkov map to show the weak coupling effect. b) The return map of one of the low-degree nodes shows that these nodes oscillate almost with the isolated local dynamics.}
    \label{return_maps_rulkov}
\end{minipage}

\subsection{H\'enon map}
\subsubsection{$u$-component coupling}
We use H\'enon map as $\bm f$ in Eq.~\eqref{laplacian-form}, a discrete-time dynamical system exhibit chaotic behavior:
\begin{eqnarray}
u_i(t+1) &=& 1 - \alpha u_i(t)^2 + v_i(t)  - \sum_{j=1}^{n} L_{ij} u_j(t)\\ 
v_i(t+1) &=& \beta u_i(t) \nonumber
\end{eqnarray}
where $\alpha=1.4$ and $\beta=0.3$ are fixed during the simulations. H\'enon maps are coupled through their $u$-components weakly. Initial positions are taken from the interval $[0,0.1]$ randomly uniformly for each map. We iterate over $11000$ steps and discard the first $10000$ steps as a transient. Our procedure reveals the local dynamics as \emph{H\'enon map}, interaction function as \emph{diffusive coupling} for $u-$component, and finally, the correct Laplacian matrix from time series data. 

\begin{figure}[H]
    \centering
    \includegraphics[width=.7\linewidth]{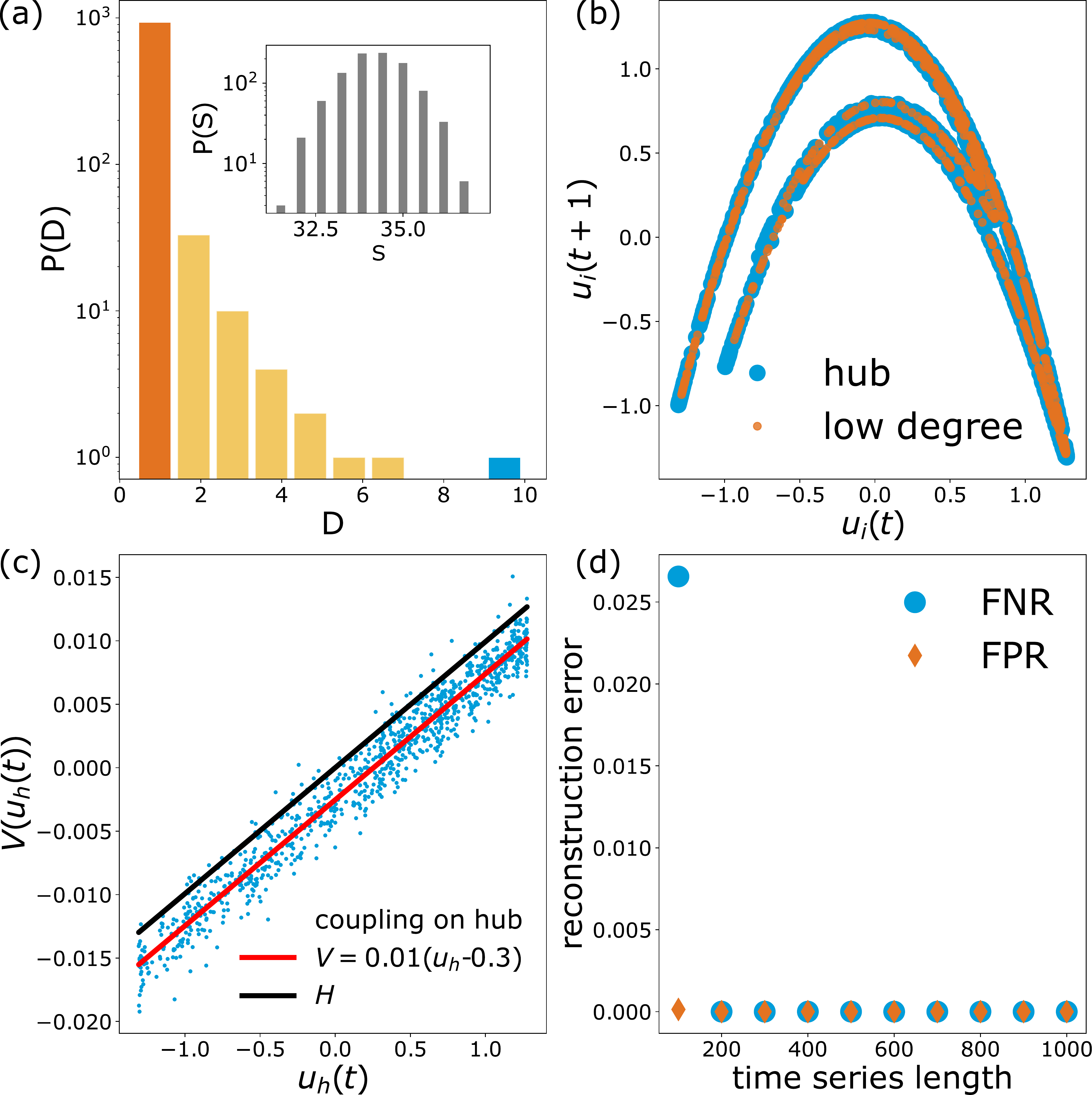}
    \caption{Reconstruction procedure for weakly $u$-coupled H\'enon maps on a real network. First, we learn a model for each node's time series data to classify them. (a) The histogram $P(D)$ based on the pairwise Euclidean distance matrix. While the highest bin of the histogram gives us $\bm f$, the lowest bin identifies the hub. The inset histogram emphasizes that the pairwise correlations of the time series observations do not provide any information about the network connectivity due to the chaotic dynamics. (b) The return maps of a low degree node and hub show the coupling effect slightly due to weak coupling. In (c), we plot the dominant coupling effect on the hub to see the model of $\bm H$. (d) shows $FNR$ and $FPR$ for different lengths of time series. $200$ data points are enough for full reconstruction of the network dynamics with $987$ nodes.}
    \label{henon-identity}
\end{figure}

\subsubsection{$u$-component coupling with sinusoidal function}

We use the same $\bm f$ but different coupling function $\bm H$ in this example: 
\begin{eqnarray}
u_i(t+1) &=& 1 - \alpha u_i(t)^2 + v_i(t)  - \sum_{j=1}^{n} L_{ij} \sin(2\pi u_j(t)) \\ 
v_i(t+1) &=& \beta u_i(t) \nonumber
\end{eqnarray}

\begin{figure}[H]
    \centering
    \includegraphics[width=.7\linewidth]{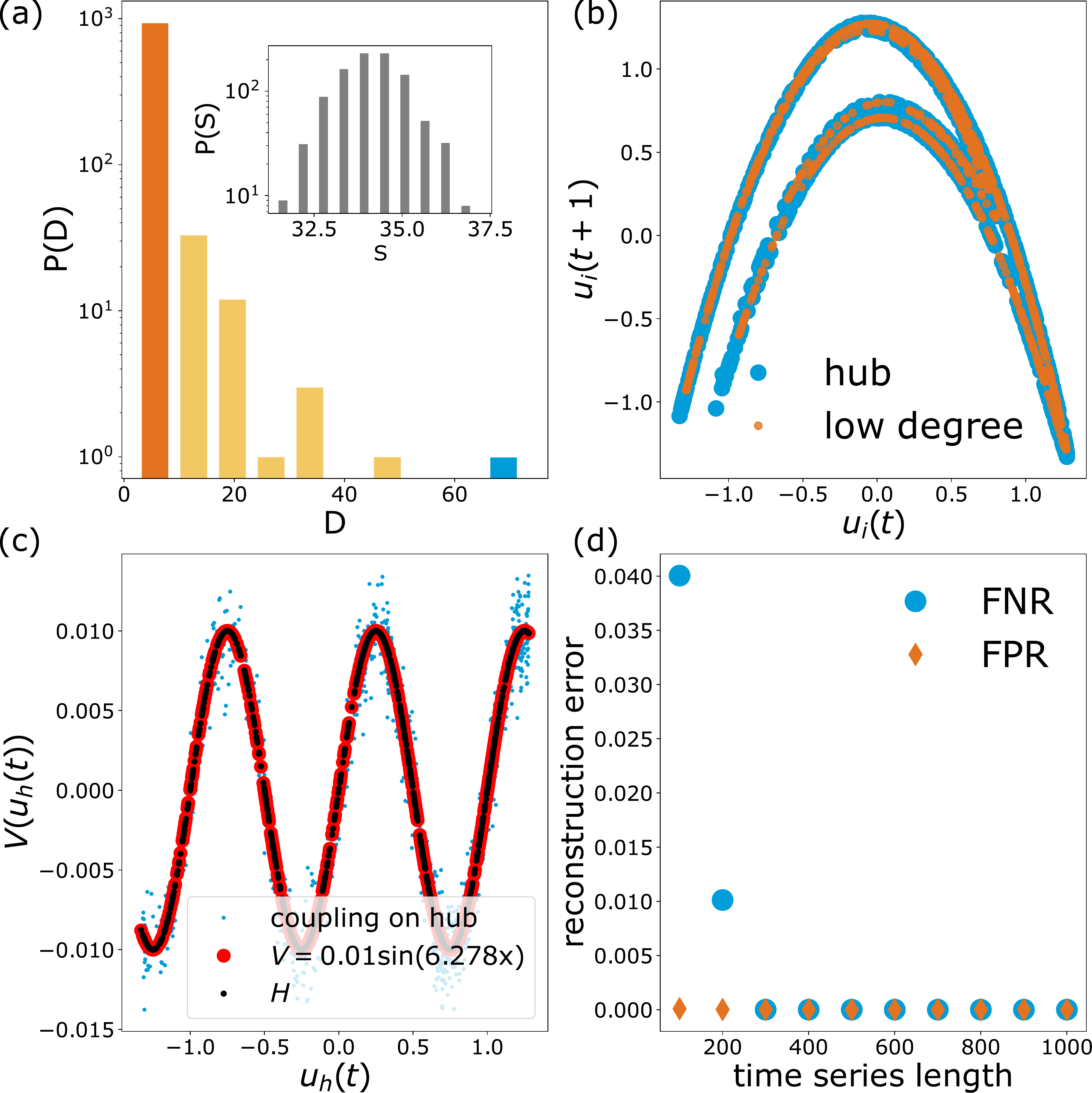}
    \caption{Reconstruction procedure for weakly $\sin{u}$-coupled H\'enon maps on a real network. First, we learn a model for each node's time series data to classify them. (a) The histogram $P(D)$ based on the pairwise Euclidean distance matrix. While the highest bin of the histogram gives us $\bm f$, the lowest bin identifies the hub. The inset histogram emphasizes that the pairwise correlations of the time series observations do not provide any information about the network connectivity due to the chaotic dynamics. (b) The return maps of a low degree node and hub show the coupling effect slightly due to weak coupling. In (c), we plot the the dominant coupling effect on the hub to see the model of $\bm H$. (d) shows $FNR$ and $FPR$ for different lengths of time series. $300$ data points are enough for full reconstruction of the network dynamics with $987$ nodes.}
    \label{henon-sine}
\end{figure}

Initial positions are taken from the interval $[0,0.01]$ randomly uniformly for each map. We iterate over $12000$ steps and discard the first $10000$ steps as a transient. Fig.~\ref{henon-sine} shows the reconstruction results of sine-coupled H\'enon maps. We simply added trigonometric functions to the candidate function library in the first step of our procedure and fully reconstructed the network connections by learning the local dynamics as \emph{H\'enon map} and the coupling function as \emph{sine function}.

\subsection{Tinkerbell Map}
Tinkerbell map coupled through their $u$-component diffusively is given by:
\begin{eqnarray}
u_i(t+1) &=& u_i(t)^2 - v_i(t)^2 + a u_i(t) +b v_i(t) - \sum_{j=1}^{n} L_{ij} u_j(t) \\ 
v_i(t+1) &=& 2u_i(t)v_i(t) + c u_i(t) + d v_i(t) \nonumber
\end{eqnarray}
where $a=0.9, b=-0.6013, c=2.0$ and $d=0.5$ are fixed during the simulations to ensure a fully chaotic regime. Data is generated from $u$-coupled Tinkerbell maps by iterating over $11000$ steps and dropping the first $10000$ time steps as a transient. We successfully reveal the local dynamics of the Tinkerbell map, interaction function and finally the Laplacian matrix from time series. 

\begin{figure}[H]
    \centering
    \includegraphics[width=.7\linewidth]{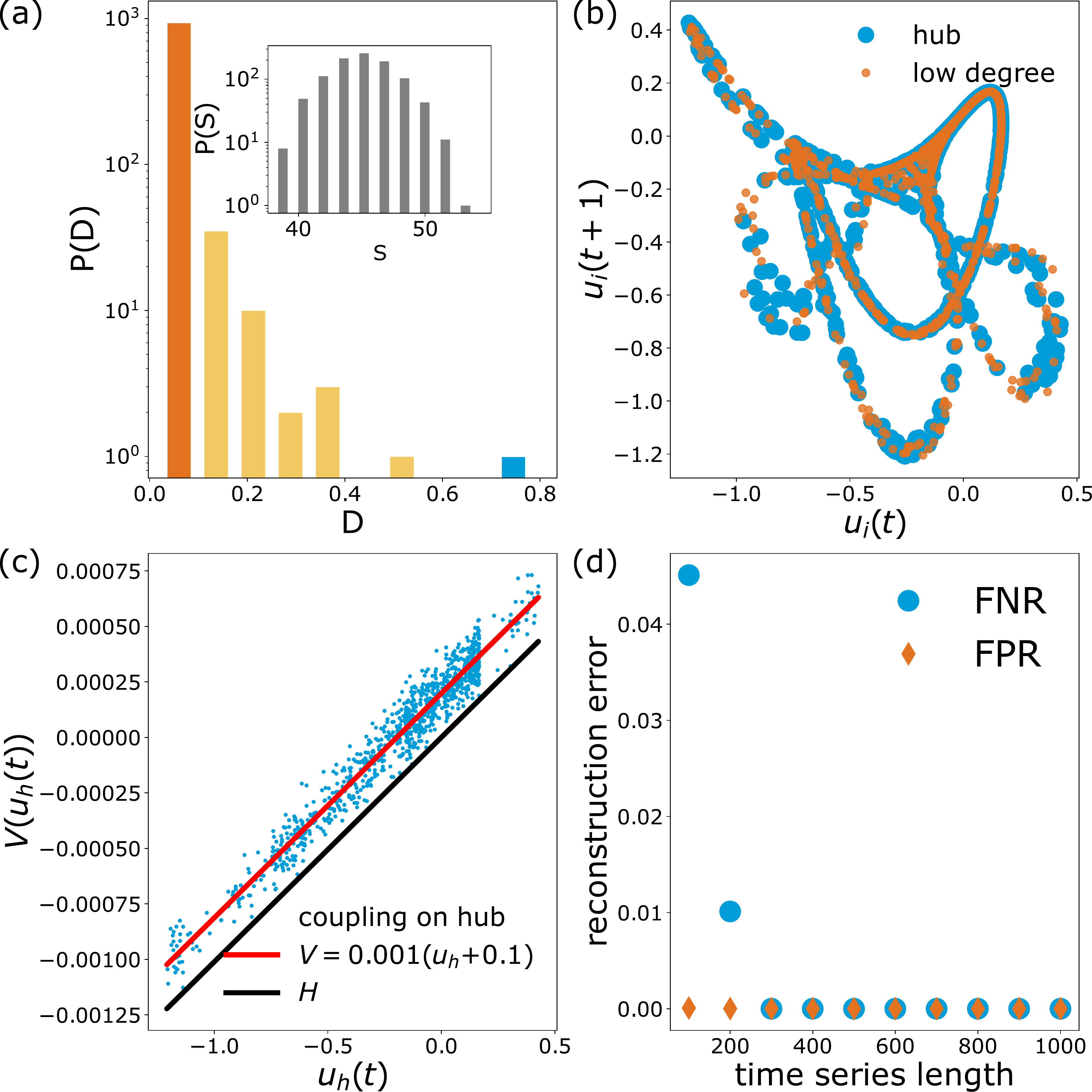}
    \caption{Reconstruction procedure for weakly $u$-coupled Tinkerbell maps on a real network. First, we learn a model for each node's time series data to classify them. (a) The histogram $P(D)$ based on the pairwise Euclidean distance matrix. While the highest bin of the histogram gives us $\bm f$, the lowest bin identifies the hub. The inset histogram emphasizes that the pairwise correlations of the time series observations do not provide any information about the network connectivity due to the chaotic dynamics. (b) The return maps of a low degree node and hub show the coupling effect slightly due to weak coupling. In (c), we plot the dominant coupling effect on the hub to see the model of $\bm H$. (d) shows $FNR$ and $FPR$ for different lengths of time series. $300$ data points are enough for full reconstruction of the network dynamics with $987$ nodes.}
    \label{tinkerbell-identity}
\end{figure}

As seen in figures \ref{henon-identity}(c), \ref{henon-sine}(c) and \ref{tinkerbell-identity}(c), the horizontal shift between the fitted curve on the data (red) and the predicted $H$ curve (black) varies for each dynamical system. This observation is compatible with the reduction theorem, which is analytically proven only for expanding maps \cite{Eroglu2020}. We show that even if the theorem is not rigorously proven for the Rulkov map, H\'enon map or Tinkerbell map, the methodology accurately works. For the coupling, $ H(x, y) = y - x$, $v(x) = \int H(x,y) dm(y) $ is found as $-x$ plus a constant. These constants are found as $1$ for coupled Rulkov maps, $-0.3$  for coupled H\'enon maps and $0.1$ for coupled Tinkerbell maps. For sinusoidal coupling, $ H(x, y) = \sin{2\pi y} - \sin{2\pi x}$, $v(x) = \int H(x,y) dm(y) $ is found as $-\sin{2\pi x}$ plus a constant and it is $0$ for $\sin$-coupled H\'enon maps.

\section{Macaque monkey visual cortex network}
\begin{figure}[H]
    \centering
    \includegraphics[width=0.8\linewidth]{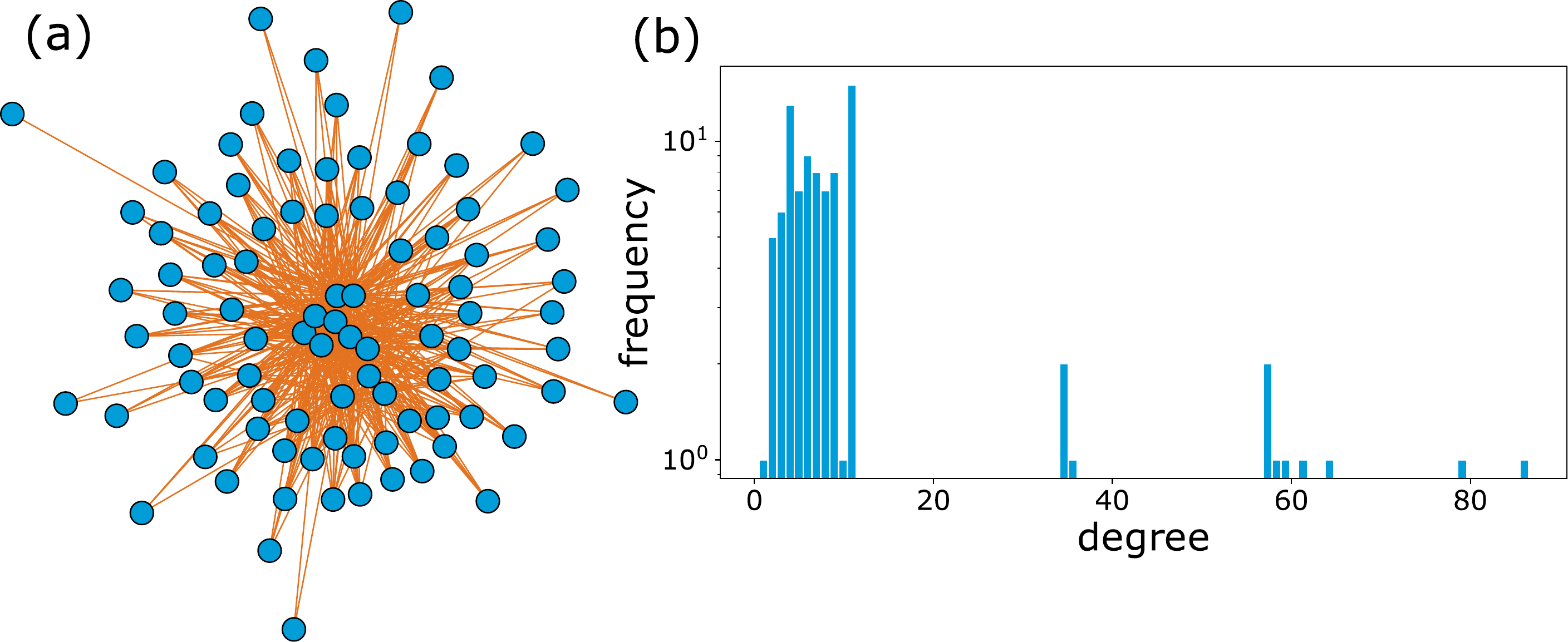}
    \caption{(a) Undirected and unweighted macaque monkey visual-cortex network. (b) Degree distribution of the associated real network.}
\label{rhesus}
\end{figure}

We tried our procedure on a real monkey visual cortex network that is not a scale-free type \cite{rhesus,Vogelstein2018}. The network consists of $91$ nodes and $581$ undirected and unweighted edges, which is illustrated in Fig.~\ref{rhesus}(a) with its degree distribution in Fig.~\ref{rhesus}(b). It is impossible to classify nodes by degree in a network with such a degree distribution, failure to classify nodes results in an unable to learn local dynamics and identify the hub. So, we assume we know the local dynamics and which node is the hub. When we continue with this preliminary information, the reduction theorem allows us to learn the coupling function smoothly and see full reconstruction.

\section{Sparse optimization}
\label{sec:regression}
The model for the dynamical system is given by the function $\bm g$:
\begin{equation}
    \bm g(\bm x) \approx \sum_{k=1}^p \psi_k(\bm x) \xi^k
\end{equation}
where $\bm \Psi(\bm x) = [\psi_1(\bm x), \psi_2(\bm x),...,\psi_p(\bm x)]$ is the library of basis functions. Only a few active terms characterize data that comes from natural systems on a well-chosen basis; most of the coefficients $\xi^k$ are zero. To find interpretable models sparsity concept is beneficial. This idea was used in a compressed sensing framework for the first time as the sparsity requirement is satisfied \cite{Wang2011-prl}. If we have prior knowledge to decide the basis, sparse regression is instrumental in avoiding overfitting and noise \cite{Brunton2016}. In our reconstruction procedure, we use sparse regression in more than one step.\\

The least absolute shrinkage and selection operator (LASSO) \cite{Tibshirani1996} is a sparse regression method that uses $\ell_1$-norm to promote sparsity. Sequentially-thresholded least-squares (STLS) is another method implemented in the Python package of SINDy \cite{desilva2020pysindy}. In STLS, we start with an ordinary least squares solution which results in overfitting the time series at each point. We repeat least squares to the values obtained by thresholding the residues below a particular cut-off value called sparsity parameter at each sequence. Finally, an interpretable dynamical model is obtained.

\section{Library of basis functions}
\label{sec:library}

Although prior knowledge about the dynamical system is always helpful in constructing a basis library, it is impossible to access the information for all cases. Nevertheless, many functions can be accurately estimated using high-order polynomials. For instance, in this work, we mainly studied with the Rulkov maps defined as follows,
\begin{eqnarray}
u(t+1) &=& \frac{\beta}{1+u(t)^2} + v(t) \label{isolated_rulkov} \\ 
v(t+1) &=& v(t) - \nu u(t) - \sigma \nonumber
\end{eqnarray}

\begin{figure}[H]
    \centering
    \includegraphics[width=\linewidth]{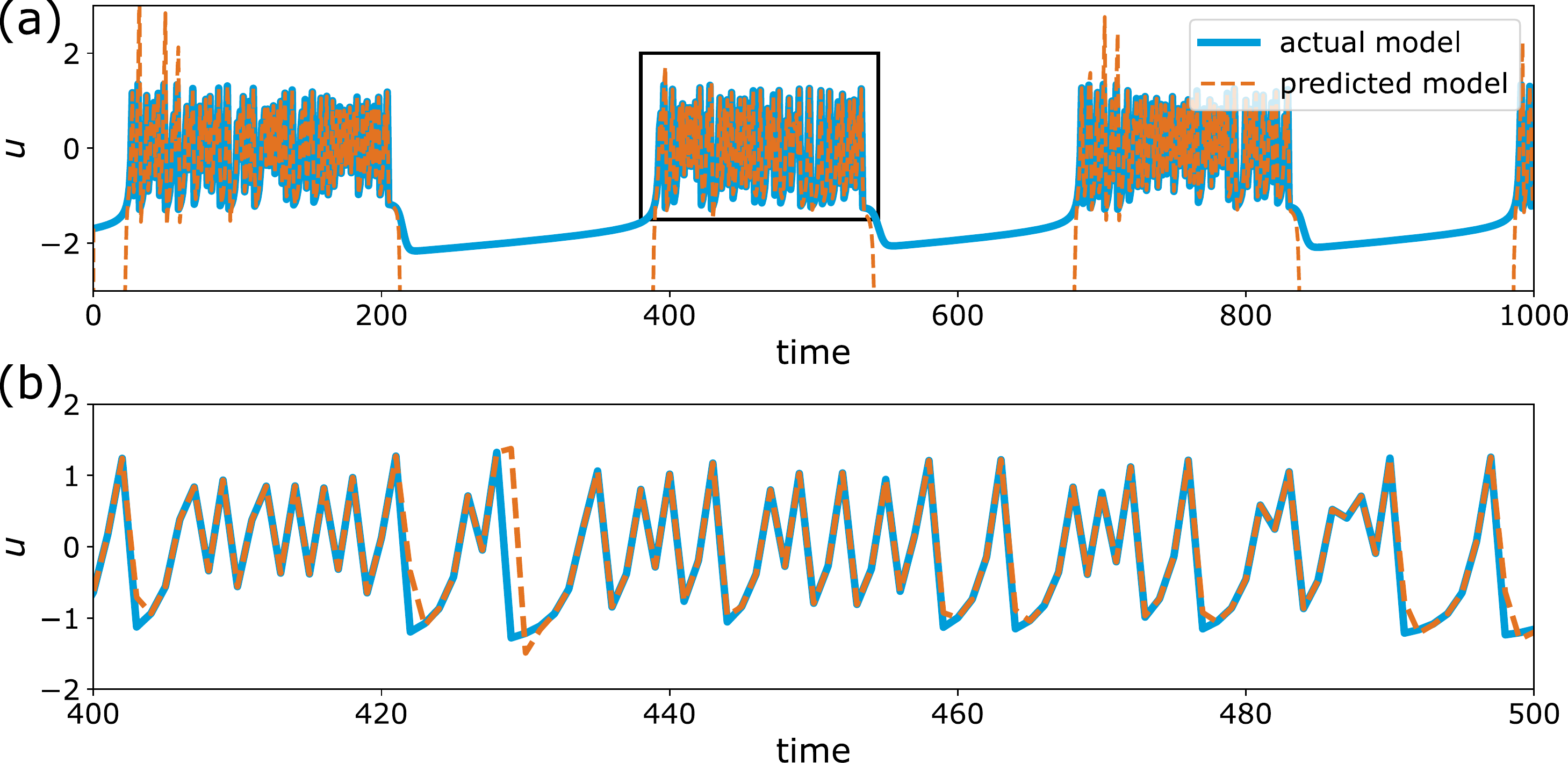}
    \caption{The predicted time series of a reconstructed Rulkov model using a candidate functions library containing only polynomials. a) predicted $u$-variable of the model b) a bursting regime package zoomed from the black rectangle in a).}
    \label{rulkov-poly}
\end{figure}

Using a basis library containing only polynomials, it is possible to reconstruct an imperfect  model for the Rulkov map, which generates quite an accurate time series to the original one (Fig.~\ref{rulkov-poly}). For the Rulkov map, only the non-polynomial term is the rational one, $\frac{1}{1+u^2}$, which is absent in the library. However, this rational term  can be expanded in a power series as follows, 
\begin{equation}
    \frac{1}{1+x^2} = \sum_{n=0}^{\infty} (-x^2)^n = 1 - x^2 + x^4 - x^6 + x^8 + \cdots 
\end{equation}
therefore, it is possible to model the Rulkov map with polynomials roughly. As the power expansion's convergent interval is $[-1,1]$, the predicted model diverges when the orbit is out of the given interval. However, such convergence problems can be solved by normalizing the time series into the associated space. The interesting dynamics for our study is the bursting regime which can be modeled without the normalization Fig.~\ref{rulkov-poly}~(b); therefore, we skip this step, but the model becomes very complicated with many higher-order polynomial terms.
In order to avoid complicated models containing many polynomial terms and possible inaccurate predictions, we prepared a standard library including common functions, such as geometric series and trigonometric functions. Although unreasonably increasing the number of functions in the library is  undesired, adding such nonpolynomial terms generally terminates a large number of polynomials. Thus, our computational cost mostly becomes cheaper using the following standard candidate library:
$$
\bm \Psi = \bigg[\bm 1, \mathcal{P}_{d}(\bm x), \big\{ \sin(i \bm x), \cos(i \bm x) \big\}_{i=1}^{r}, \Big\{ \frac{1}{\bm{x}^i}, \frac{1}{1 \pm \bm x^i}, \frac{1}{(1 \pm \bm{x})^i} \Big\}_{i=1}^{q} \bigg]
$$
where $d$ is the degree of the polynomials, and for instance $\mathcal{P}_{2}(\bm x)$ is:
$$
\mathcal{P}_2 (\bm x) = [u, v, u^2, uv, v^2]. 
$$
Our main strategy is to start $d=r=q=1$, and increase them slowly to estimate the model with the minimum number of candidate functions. The first step of our reconstruction algorithm aims to learn the local dynamics by classifying the inferred model equations due to their similarity. Here, we illustrate the reconstructed models on the mouse neocortex network using coupled Rulkov maps (Table~\ref{inferred_equations}) through a set of selected nodes' with respect to their degrees $k_i$ ($k_0$ = 0, $k_{46}$ = $k_{136}$ = 4, and $k_{218}$ = 26). As node-$0$ has no incoming links, we reconstructed only the Rulkov map dynamics. Node-$46$ has $4$ incoming links, the inferred coefficients for it do not reflect the exact governing equation, but one can say that its governing equations are similar to node-$136$, which also has $4$ incoming links. The most \emph{distinct} model belongs to node-$218$, which is the hub of this particular network.

\begin{minipage}[H]{\linewidth}
\centering
    \begin{tabular}{|p{1cm}|p{1cm} p{1cm} p{1cm} p{1cm} p{1cm} p{1cm} p{1cm}|}
    \hline
    node id & $1$ & $u$ & $v$ & $\cdots$ & $\cos(u)$ & $\frac{1}{1-u}$ & $\frac{1}{1+u^2}$ \\
    \hline
     0 & -0.000 &  -0.000  & 1.000 & \vdots & 0.000  & -0.000 &  4.100 \\
    \vdots & \vdots & \vdots & \vdots & \vdots & \vdots  & \vdots & \vdots  \\
    46 & 0.000 &  -0.015 & 1.005& \vdots & 0.005  & 0.000 &  4.093 \\
    \vdots & \vdots & \vdots & \vdots & \vdots & \vdots &  \vdots & \vdots \\
    136 & -0.000 &  -0.015  & 1.006 & \vdots & 0.000  & -0.000 &  4.101 \\
    \vdots & \vdots & \vdots & \vdots & \vdots  & \vdots & \vdots  & \vdots \\
    218 & -0.714 &  0.098 & 0.783 & \vdots & 0.016  & 0.001 &  4.075 \\
    \hline
    \end{tabular}
    \captionof{table}{Table for the inferred coefficients of the corresponding candidate functions for some of the selected nodes. We present only nonzero functions for any node.}
    \label{inferred_equations}
\end{minipage}

The first step of the reconstruction algorithm is to classify nodes, which uses a distance between nodes' coefficients. We define the distance metric as,
$d_{ij}=(\Sigma_{k=1}^{p}~\frac{1}{V_k}~|\xi_i^k~-~\xi_j^k|^2 )^{1/2}$ 
where $|\cdot|$ is absolute value, $V_k$ is the variance of the predicted coefficients of the $k$-th function in $\bm \Psi$. As an example, we compute the pairwise distances of the learnt models given in Table~\ref{inferred_equations}:

\begin{minipage}[H]{\linewidth}
\centering
    \begin{tabular}{l|rrrr}
    node id &  0   &  46   &  136 &  218 \\
    \hline
    0   &   0.000 &   2.830 &  2.105 &  18.351 \\
    46   &   2.830 &   0.000 &  1.925 &   16.808 \\
    136 &  2.105 &  1.925 &   0.000 &  16.979 \\
    218 &   18.351 &   16.808 &  16.979 &   0.000 \\
    \end{tabular}
    \captionof{table}{Pairwise distances of the predicted models for some of the selected nodes.}
    \label{distance_matrix}
\end{minipage}

As numerically shown in Table~\ref{distance_matrix}, $d_{ij}$ is small for the nodes with the same or similar degrees, such as the distance is 1.925 between $46~(k_{46}=4)$ and $136~(k_{136}=4)$, while the distance is large for different predicted models, such as the distance is 18.351 between $46~(k_{46}=4)$ and $218~(k_{218}=26)$.  

\section{Comparison between sparse regression methods}
We compare two optimizers to perform sparse regression; LASSO and STLS.\\

LASSO is defined by the equation \cite{Mehta2019a}:
\begin{eqnarray}
\bm{\hat w}_{LASSO}(\lambda) = \underset{\bm w \in \mathbb{R}^{m}}{\mathrm{argmin}} (||\bm X \bm w - \bm y||_2^2 + \lambda ||\bm w||_1)
\end{eqnarray}
for a linear problem $\bm y = \bm X \bm w$ where $\bm y \in \mathbb{R}^n$, $\bm X \in \mathbb{R}^{n \times m}$ and $\lambda$ is \emph{penalty term} and uses $\ell_1$-norm to penalize the weights.\\

In STLS, Ridge regression is used and defined by the equation \cite{Mehta2019a}:
\begin{eqnarray}
\bm{\hat w}_{Ridge}(\lambda) = \underset{\bm w \in \mathbb{R}^{m}}{\mathrm{argmin}} (||\bm X \bm w - \bm y||_2^2 + \lambda ||\bm w||_2^2)
\end{eqnarray}
and $\lambda$ is fixed as $0.05$. Sparsity is obtained by masking out elements of $\bm {\hat w}$ that are below a given threshold. This threshold is called as \emph{sparsity parameter} and it is the single hyper-parameter of STLS. In the final step of our reconstruction approach, we find a sparse solution for $\bm w$, it corresponds to $\bm L$ in our case. 
From the algorithmic point of view, fine-tuning of the hyper-parameters, penalty term for LASSO and sparsity parameter for STLS, is essential to improve the accuracy. These two parameters are not equivalent, however we can make a similar interpretation between them. The length of the time series determines the problem type since our network size is fixed in the experiments. The time series length lower than the network size $(T < nm = 987 \times 2=1974)$ correspond to an under-determined case for the real network we use.\\ 

\begin{figure}[H]
    \centering
    \includegraphics[width=.9\linewidth]{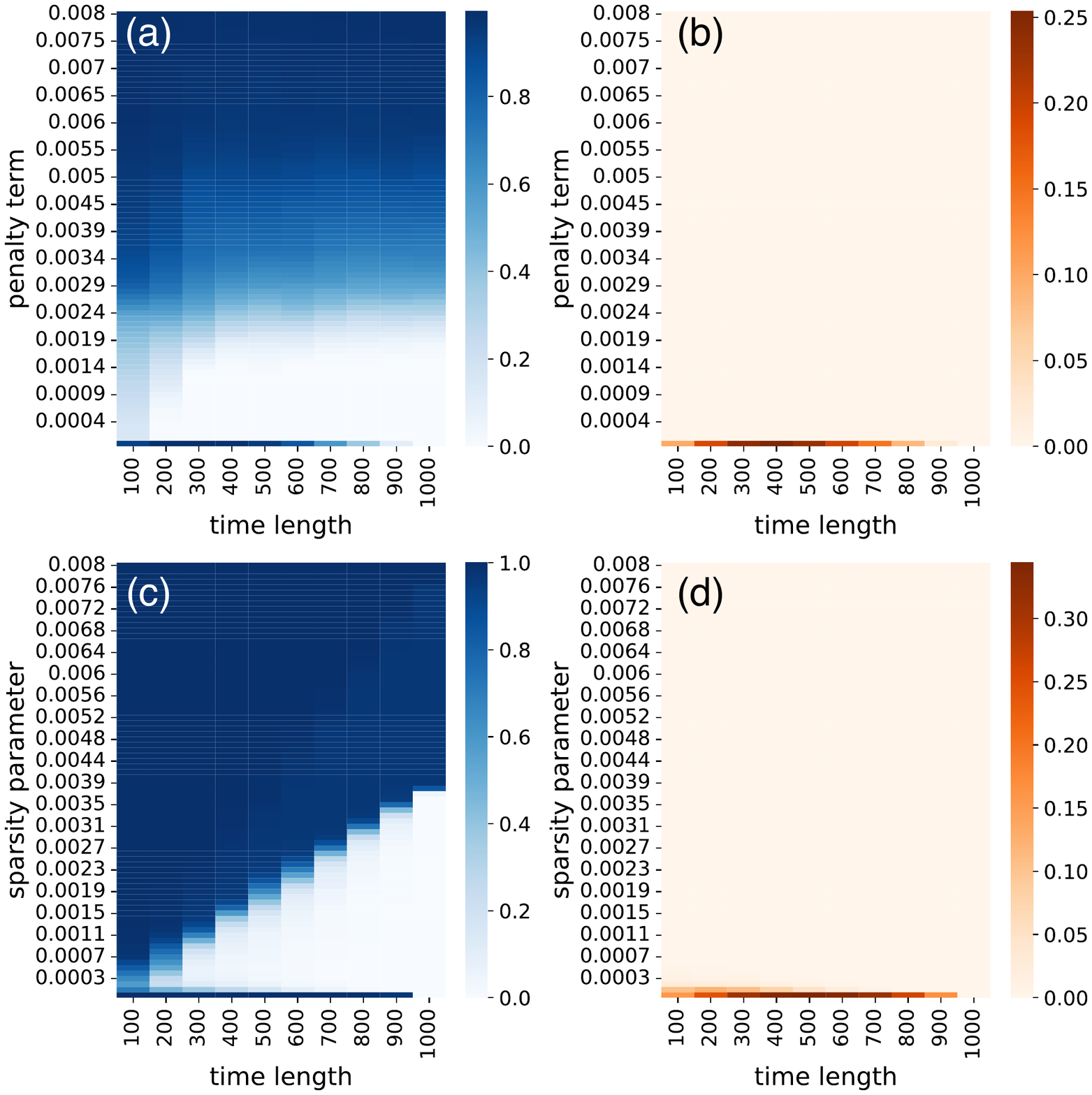}
    \caption{(a) $FNR$ and (b) $FPR$ with respect to a list of penalty terms for LASSO. (c) $FNR$ and (d) $FPR$ with respect to a list of sparsity parameters for STLS. $10$ different lengths of time series are used to compare the performance.}
    \label{lasso-stls-heatmap}
\end{figure}

\begin{figure}[H]
    \centering
    \includegraphics[width=0.8\linewidth]{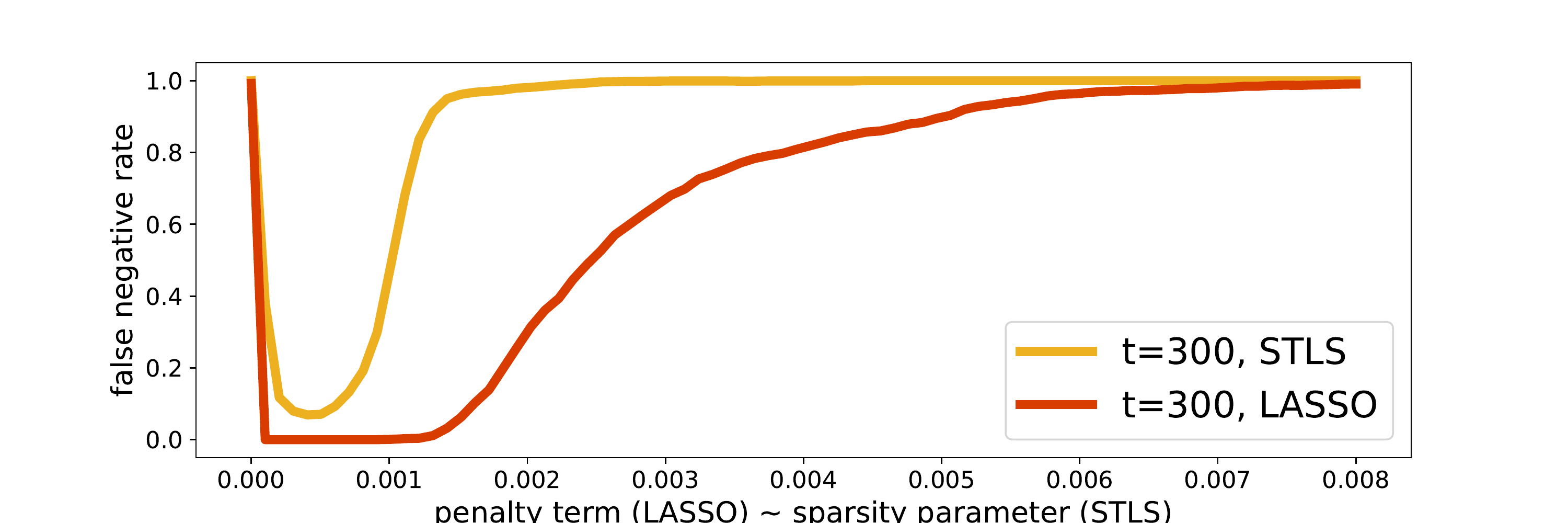}
    \caption{$FNR$ with respect to a list of penalty term for Lasso and sparsity parameter for STLS. $(T = 300)$ case is used as a short time series example. There is no \emph{full reconstruction} region for STLS.}
    \label{lasso-stls-example}
\end{figure}

Fig.~\ref{lasso-stls-heatmap} presents the reconstruction performance for different time series lengths and two hyperparameters as $FNR$ and $FPR$. We compare $FNR$ values to evaluate the performance of LASSO and STLS in Fig.~\ref{lasso-stls-heatmap}(a) and Fig.~\ref{lasso-stls-heatmap}(c), since $FPR$ values are always very small. The results show that LASSO overcomes STLS for shorter time series. Since the $\ell_1$-norm penalized solution, is a generalization of compressive sensing approach, it provides a unique solution for underdetermined cases \cite{donoho}. On the other hand, if there is long enough data, STLS becomes a fast alternative to the LASSO from an algorithmic point of view and works very well for our problem. We took one under-determined case $(T = 300)$ to compare the performance of two different sparsity promoting methods in Fig.~\ref{lasso-stls-example}. If we have limited data, we see full reconstruction by LASSO but not by STLS. \\


\section{Effect of network sparsity on reconstruction performance}

\begin{figure}[H]
    \centering
    \includegraphics[width=\linewidth]{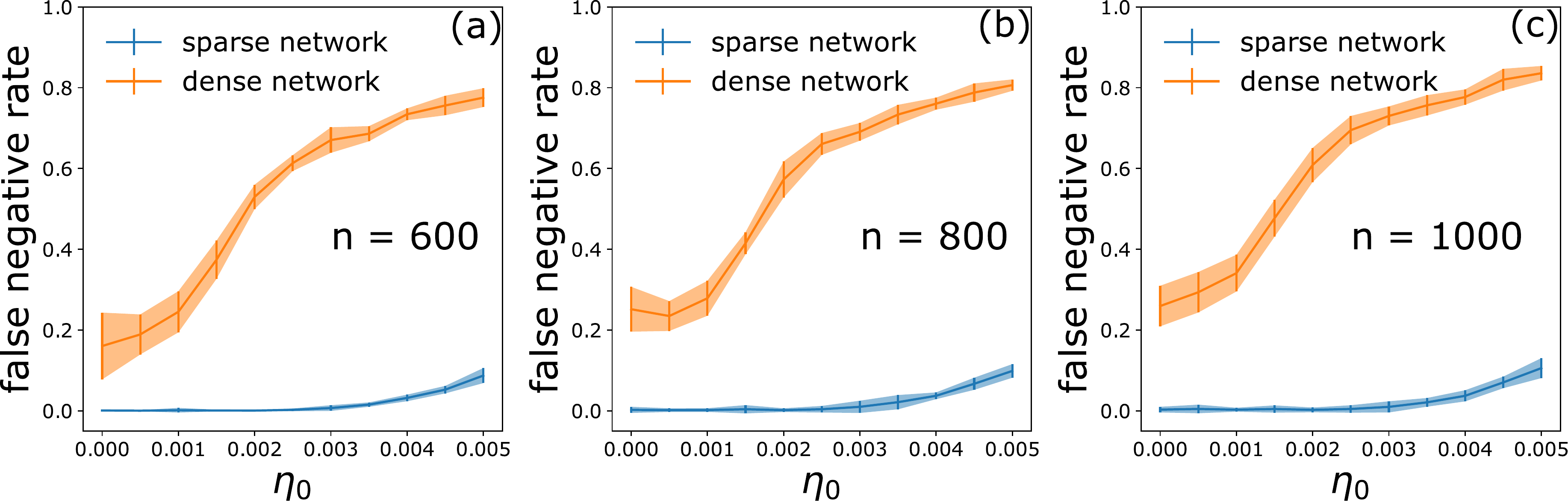}
    \caption{$FNR$ against noise for sparse and dense networks of sizes (a) $n=600$, (b) $n=800$ and (c) $n=1000$. All points show an average error over $10$ different realizations of scale free networks and shaded regions present standard deviations. Time series length is $500$, which corresponds to underdetermined cases for all three network sizes. }
    \label{sparse-dense-network}
\end{figure}

In addition to the lack of sufficient data, the sparsity condition of the Laplacian matrix should be satisfied to solve the linear regression problem within the compressed sensing framework. We present results in this direction based on synthetical networks. The algorithm we use to generate directed scale-free networks has three probabilistic parameters: $\alpha$ and $\gamma$ are the probabilities of adding a new node connected to an existing node chosen randomly according to the in-degree and out-degree distribution, respectively. $\beta$ is the probability of adding an edge between two existing nodes. We use two different parameter setting as $[\alpha = 0.41, \gamma = 0.05, \beta = 0.54]$ and $[\alpha = 0.2, \gamma = 0.5, \beta = 0.3]$  to obtain Laplacian matrices at different sparsity levels. The former (default) setting's power exponent is approximately $2$ and generates denser networks, while the latter's power exponent is approximately $1$ and generates sparser networks. As seen in Fig.~\ref{sparse-dense-network}, we see full reconstruction in sparser networks since the Laplacian matrix meets the sparsity condition for the compressive sensing approach \cite{candes2006stable}.

\section{Cross validation of inferred Laplacian matrices}
In the main text, we use the ground-truth Laplacian matrix to show the performance of the reconstruction for various hyper-parameters and data lengths. However, it is unlikely to know the true matrix to evaluate the methodology's success in real-world problems. Therefore, it is important to assess the reconstruction outcome when the real matrix is absent using a cross-validation technique. To illustrate the cross-validation approach, we obtained \emph{predicted time series} of length 500 using inferred Laplacian matrices for various hyper-parameters (penalty term) $\lambda=0.008, 0.003$ and $0.0003$. Out of those three parameters, the outcome for $\lambda=0.0003$ gives the best-fit, Fig.~\ref{test_series}(a).

\begin{minipage}[H]{\linewidth}
    \centering
    \includegraphics[width=\linewidth]{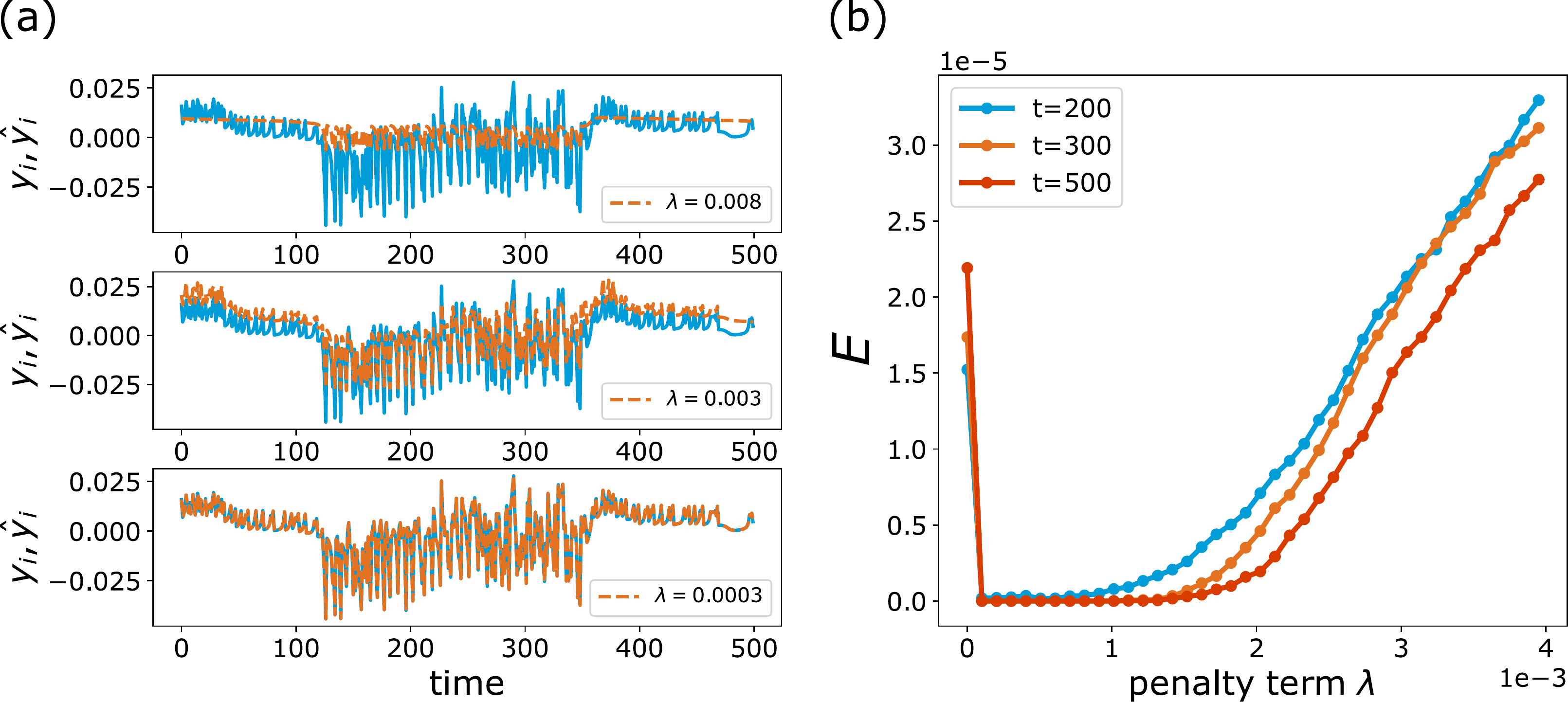}
    \captionof{figure}{a) Comparison between measurement of a node and test time series. Each $\lambda$ parameter corresponds to different connectivity matrices. We simulate the inferred models and obtain test time series. b) Cross validation of the inferred time series based on obtained connectivity matrices for a series of penalty terms.}
    \label{test_series}
\end{minipage}

To find the best hyper-parameters in large network dynamics, we use an error  function
\begin{eqnarray}
    E = \frac{1}{n}\sum_{i=1}^n \Bigg( \frac{1}{T}\sum_{t=1}^T \Big( y_{i,t} - \hat{y}_{i,t} \Big)^2 \Bigg)
\end{eqnarray}
where $y_i = x_i(t+1) - f(x_i)$ and $\hat{y}$ denotes the predicted one. Then we optimize the hyper-parameters, which minimize the error. The parameters yielding the most accurate results can be used as the optimal hyper-parameter (Fig.~\ref{test_series}(b)).

\section{Why is the reduction theorem crucial to reconstruct large networks?}
\label{sec:whyreduction}
A fully-algorithmic approach is to apply any sparse optimization technique directly to the entire multi-variable time series. Only up to a certain size networks can be recovered by this method. By the reduction theorem, we learn local dynamics from low degree nodes and interaction dynamics from the hub. We are not looking for a sparse linear combination of functions within a massive library of all possible functions. Therefore, there is no limit on network size in our approach. In this section, we will show this with a simple example:\\

Assume we have a network with $5$ nodes. If node $1$ has links coming from the nodes $0$, $4$ and $5$ as seen in Fig.~\ref{5_net}, then the dynamics of node $1$ is written simply in Laplacian matrix form:
\begin{figure}[H]
    \centering
    \includegraphics[width=0.2\linewidth]{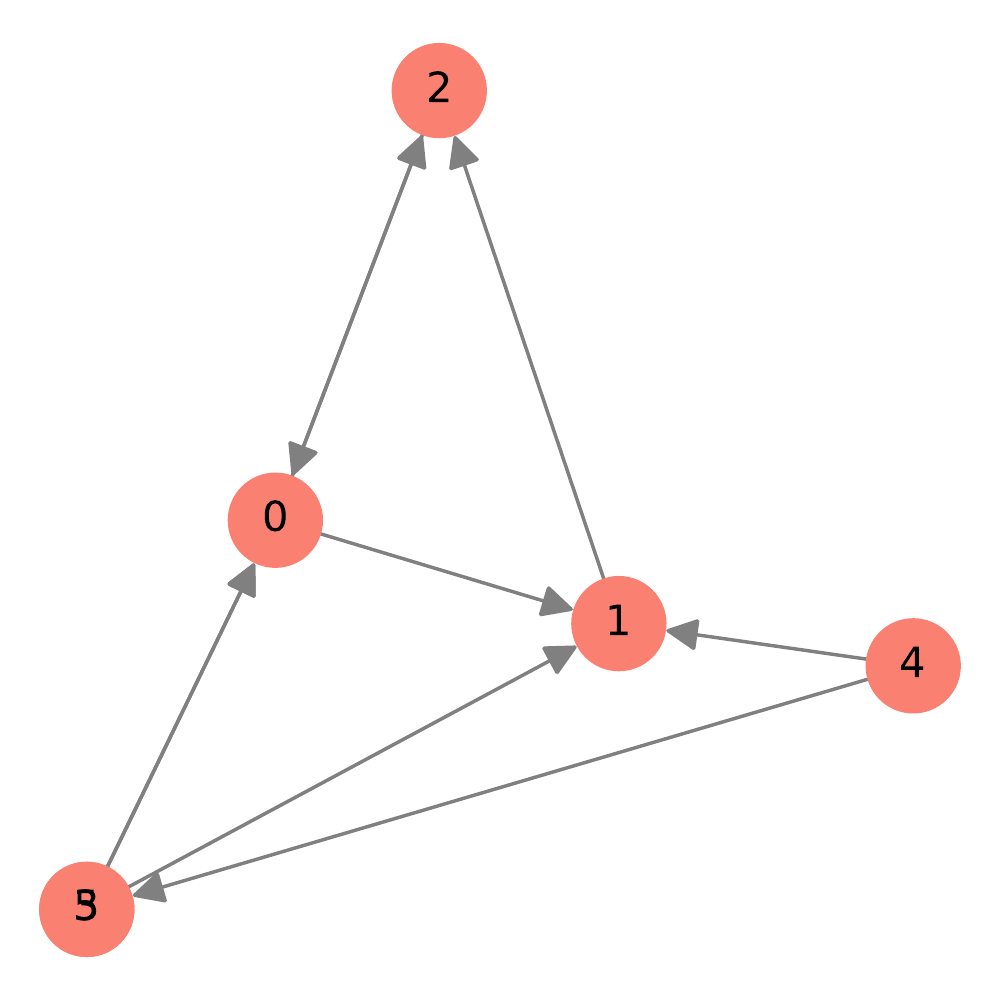}
    \caption{$5$-nodes network as an example.}
    \label{5_net}
\end{figure}
\begin{eqnarray}
\bm{x}_1(t+1) = \bm{f} (\bm{x}_1(t)) + 
3\bm{H}(\bm{x}_1(t))  - \bm{H}(\bm{x}_0(t)) - \bm{H}(\bm{x}_4(t)) - \bm{H}(\bm{x}_5(t))
\label{explicit-form-laplacian}
\end{eqnarray} 
where $\bm x_i \in \mathbb{R}^m$. Eq.\eqref{explicit-form-laplacian} is in the form of a linear combination of the connections, allowing any sparse regression algorithm to capture the parameters correctly in case of a well-defined basis library. For one single equation, it seems pretty feasible, but for a $5$-nodes network, each has $m$ components, we have $5m$ equations to be discovered. Suppose we set the basis composed of polynomials up to second-degree and all quadratic combinations of all variables due to a piece of prior knowledge about the coupling function. This example setting has $9$ functions for $m=3$. In that case, our basis includes $45$ features for this particular example. If we do not have prior knowledge, we use an extensive library that includes more candidate functions. We performed experiments in this direction to see the limitations. The results showed that a purely algorithmic equation-based learning approach fails as the network size, that is, the basis grows. Results of Novaes \emph{et al.} validated our conclusions, they showed that adding new functions to the basis due to increase of network size (up to $20$-nodes) effects badly the reconstruction \cite{Novaes2021}.\\
Furthermore, in the final sparse regression step, our problem turns into a linear equation as $\bm y = \bm A \bm x$, where $\bm y \in \mathbb{R}^T$, $\bm A \in \mathbb{R}^{T \times n}$ and $\bm x \in \mathbb{R}^n$. Here, $\bm y$ is the data, and the length of the data is $T$, $n$ is the network size, and we focus on the underdetermined cases ($T < n$) since data availability is always limited.
From the perspective of compressed sensing, the measurement matrix $\bm A$ should satisfy some conditions to converge the sparsest solution for $\bm x$ by $\ell_1$-norm. One of these conditions can be summarized as the time series length scales with the network size, so if the measurement matrix has more columns, longer measurements are needed. It is given by the relation:
$$
T \approx \mathcal{O} (K log(n/K))
$$
where $K$ is the network's sparsity measure (number of nonzero entries) \cite{brunton_kutz_2019}. If the number of columns is determined by the basis functions that include all possible pairwise interactions, not the network size, this means the need for a longer time series.
On the other hand, the uncertainty of the interaction function is an obstacle to inferring the exact connection matrix \cite{Liu2016}. Finding $\bm H$ and therefore reducing the dimension of the optimization problem for full reconstruction broadens the applicability and distinguishes our work.

\section{General Diffusive Coupling and Master stability function}
\label{sec:msf}
Consider the function
$\bm h: \mathbb{R}^m \times \mathbb{R}^m \rightarrow \mathbb{R}^m$. We say that $\bm h$ is diffusive if  

$$
\bm h(x,x) = \bm h(0) = 0 \mbox{      and         } \bm h(x,y) = - \bm h(y,x)
$$

We again consider the model to a general diffusive coupling

\begin{equation}
\bm{x}_i(t+1) =\bm{f}(\bm{x}_i(t)) + \gamma\sum_{j=1}^n w_{ij} \bm h (\bm{x}_j,\bm{x}_i)
\label{eq:motion_lap}
\end{equation}
where $\gamma$ is the coupling strength.

We perform the analysis close to synchronization $\bm{x}_i = \bm s + \bm{\xi}_i$ so
$$
\bm h(\bm{x}_j,\bm{x}_i) = \bm h(\bm s+ \bm{\xi}_j, \bm s + \bm{\xi}_i) = \bm h(\bm s,\bm s) + D_1 \bm h(\bm s, \bm s) \bm{\xi}_j + D_2 \bm h(\bm s, \bm s) \bm{\xi}_i 
$$
but since the coupling is diffusive
$$
D_2 \bm h(\bm s,\bm s) = - D_1 \bm h(\bm s,\bm s)
$$
we get
$$
\bm h(\bm{x}_j,\bm{x}_i) = \bm{H}(\bm s)(\bm{\xi}_j - \bm{\xi}_i) + \bm R(\bm{\xi}_i,\bm{\xi}_j)
$$
where $\bm H(\bm s)  = D_1 \bm h(\bm s, \bm s)$ and  $\bm R(\bm{\xi}_i,\bm{\xi}_j)$ contains quadratic terms. Then, the first variational equation about the synchronization manifold
\begin{eqnarray}
\bm{\xi}_i(t+1) &=& D\bm{f}(\bm{s}(t)) \bm{\xi}_i(t) + \gamma \sum_{j=1}^n w_{ij} \bm H(\bm s) (\bm{\xi}_j(t) - \bm{\xi}_i(t)) \\
&=&D\bm{f}(\bm{s}(t)) \bm{\xi}_i(t) - \gamma \bm H(\bm s) \sum_{j=1}^n L_{ij} \bm{\xi}_j(t). 
\label{Dif}
\end{eqnarray}
All blocks have the same form which are different only by $\lambda_i$, the $i$th eigenvalue of $\bm L$. Then we obtain the parametric equation for the modes
$$
\bm{u}(t+1) = [D\bm f(\bm s(t)) - \kappa \bm H(\bm s(t))]\bm{u}(t) 
$$ 
where $\kappa = \gamma \lambda_i$. By fixing $\kappa$, we compute the maximum Lyapunov exponent  $\Lambda(\kappa)$ as 
$$
\| \bm u(t) \| \le C e^{\Lambda(\kappa) t}.
$$
The map 
\begin{equation}\label{msf}
\kappa \mapsto \Lambda(\kappa)
\end{equation}
is called master stability function. Notice that if  $\Lambda (\kappa) <0$ when $\kappa > \gamma_c \lambda_2$ then 
$\|u\| \rightarrow 0$. 
For more details on the master stability function Eq. \ref{msf}, see e.g., \cite{PhysRevE.80.036204, 2016Phs}.

\bibliography{paper}

\end{document}